\magnification=\magstep1
\input amstex
\voffset=-3pc
\documentstyle{amsppt}
\NoBlackBoxes
\vsize=9truein
\hsize=6.5truein

\def\script U{\Cal U}
\def\script W{\Cal W}
\def\script S{\Cal S}
\def\cK{\Cal K}
\def\cV{\Cal V}
\def\bR{\Bbb R}
\def\bC{\Bbb C}
\def\var{\varepsilon}
\def\her{\text{ her}}

\topmatter
\title
Nearly Relatively Compact Projections in Operator Algebras
\endtitle

\rightheadtext{Nearly Relatively Compact Projections}
\leftheadtext{Lawrence G.~Brown}
\author
Lawrence G. Brown\\
Department of Mathematics\\
Purdue University\\
W. Lafayette, IN 47907, USA\\
e-mail: lgb \@ math. purdue.edu
\endauthor
\abstract Let $A$ be a $C^*$-algebra and $A^{**}$ its enveloping von Neumann algebra.
C.~Akemann suggested a kind of non-commutative topology in which certain projections in $A^{**}$ play the role of open sets.
The adjectives ``open'', ``closed'', ``compact'', and ``relatively compact'' all can be applied to projections in $A^{**}$.
Two operator inequalities were used by Akemann in connection with compactness.
Both of these inequalities are equivalent to compactness for a closed projection in $A^{**}$, but only one is equivalent to relative compactness for a general projection.
A third operator inequality, also related to compactness, was used by the author.
It turns out that the study of all three inequalities can be unified by considering a numerical invariant which is equivalent to the distance of a projection from the set of relatively compact projections.
Since the subject concerns the relation between a projection and its closure, Tomita's concept of regularity of projections seems relevant, and some results and examples on regularity are also given.
A few related results on semicontinuity are also included.
\endabstract
\endtopmatter


\noindent
{\bf AMS(MOS) subject classifications (1985 revision)}.
Primary 46L05, Secondary 47C15

\document
\subheading{\S 1. Introduction}

A projection in $A^{**}$ is called \underbar{open} if it is the support projection of a hereditary $C^*$--subalgebra of $A$.
$p$ is \underbar{closed} if $1-p$ is open.
$Q(A)$, the quasi-state space of $A$, is $\{f\in A^*\colon f\geq 0$ and $\|f\|\leq 1\}$.
$S(A)$, the state space of $A$, is $\{f\in Q(A)\colon \|f\|=1\}$.
For a projection $p$ in $A^{**}$, let $F(p)=\{f\in Q(A)\colon f(1-p)=0\}$.
Then $p$ is closed if and only if $F(p)$ is weak* closed (Effros [13]).
$p$ is called \underbar{compact} if $F(p)\cap S(A)$ is weak* closed.
For every projection $p$ in $A^{**}$, there is a smallest closed projection $\overline p$ such that $\overline p\geq p$.
$\overline p$ is called the \underbar{closure} of $p$.
$p$ is called \underbar{relatively compact} if $\overline p$ is compact.
For any subset $S$ of $A^{**}$, $S_{sa}$ denotes $\{x\in S\colon x=x^*\}$ and $S_+$ denotes $\{x\in S\colon x\geq 0\}$.
If $A$ has a unit, then every projection in $A^{**}$ is relatively compact.
Therefore our concern is with non-unital $C^*$-algebras.

Consider the following properties for a projection $p$ in $A^{**}$:

\itemitem{(1)}$\exists a\in A_{sa}$ such that $p\leq a\leq 1$.
\itemitem{(2)}$\exists a\in A_{sa}$ such that $p\leq a$.
\itemitem{(3)}$\exists a\in A_{sa}$ such that $p\leq pap$.

Clearly (1) $\Rightarrow$ (2) $\Rightarrow$ (3), and any of the properties for $\overline p$ implies the same property for $p$.
Akemann [4] showed that for $p$ closed each of (1) and (2) is equivalent to compactness and for general $p$, (1) is equivalent to relative compactness, but for general $p$, (2) does not imply relative compactness.
We showed in [8] that for $p$ open and $A$ $\sigma$-unital, (3) is equivalent to the property that every closed subprojection of $p$ is compact.
We will show below that half of this result is true for general $A$, but unfortunately nothing in this paper ``explains'' the result.

The original goal of this work was to find all possible answers to:\ Which of (1), (2), (3) are true for $p$ and which are true for $\overline p$?
There are, in fact, six possible answers, but it is better to organize the subject differently.
If a non--zero $p$ satisfies (3), let $\alpha(p)=\inf\{\|a\|\colon a\in A_{sa}$ and $p\leq pap\}$.
Otherwise, let $\alpha(p)=\infty$.
Also let $\alpha(0)=1$.
Clearly $1\leq \alpha(p)\leq\infty$ and $\alpha(p)\leq \alpha(\overline p)$.
More generally, $p_1\leq p_2\Rightarrow\alpha (p_1)\leq\alpha (p_2)$, so that $\alpha(p)$ is some kind of measure of how large $p$ is.
(2) will be shown equivalent to ``$\alpha(p)=1$'', and hence (1) is equivalent to ``$\alpha(\overline p)=1$''.
Thus all the information for our original goal is contained in the pair $(\alpha(p),\alpha(\overline p))$.
We will give enough examples to show that every pair $(s,t)$ such that $1\leq s\leq t\leq \infty$ is $(\alpha(p),\alpha(\overline p))$ for some $p$ and $A$.

Let $RC$ be the set of relatively compact projections in $A^{**}$, $ORC$ the set of open relatively compact projections, and $CRC$ the set of compact projections (which, of course, is the same as the set of closed relatively compact projections).
Then for any projection $p$ in $A^{**}$, dist$(p,RC)=[1-\alpha(p)^{-1}]^{1/2}$, where the distance is with respect to the metric induced by the norm.
Also if $p$ is open, then dist$(p,RC)=$ dist $(p,ORC)$; and if $p$ is closed, then dist$(p,RC)=$ dist$(p,CRC)$.
Now $CRC$ is a norm closed set, because of the semicontinuity characterization of compactness, [7, 2.47(iv)].
Thus dist$(p,CRC)=0$ implies $p\in CRC$.
Neither $RC$ nor $ORC$ need be closed, since Akemann's counterexample in [4] showing (2) $\not\Rightarrow$ (1) uses an open projection.
Thus in some sense our results ``explain'' Akemann's results that (2) $\Rightarrow$ (1) for closed projections but not for general projections.

A projection $p$ in $A^{**}$ will be called \underbar{nearly relatively compact} if dist$(p,RC)<1$.
By our results proved below, this is equivalent to ``$\alpha(p)<\infty$'' or ``$p$ satisfies (3)''.
We will not define ``nearly compact''.
The reader might think this should mean ``dist$(p,CRC)<1$''; but we think a better meaning for this term would be ``closed and nearly relatively compact''.
We give some discussion of this point below, but do not consider the issue to be completely settled.

There are other natural interpretations of $\alpha(p)$ which are included, together with the main results, in Section 2, except for the examples, which are in Section 3.
Section 4 contains some results and examples on regularity of projections and its relation to the above.
Section 5 contains special results on open projections, Section 6 results on $\alpha(p_1\vee p_2)$, and Sections 7,8,9 contain miscellaneous related results, remarks, and examples.

\subheading{\S 2.\ $\alpha(p)$}

\proclaim{Theorem 2.1}If $p$ is a projection in $A^{**}$, then $\alpha(p)=1$ if and only if $p\leq a$ for some $a$ in $A_{sa}$.
\endproclaim

\demo{Proof}We rely on a result of Akemann, Theorem 1.2 of [3], which states in slightly different words:\ If $A$ is a $C^*$-subalgebra of $B$ and $c$ is a positive element of her$_B(A)$, the hereditary $C^*$-subalgebra of $B$ generated by $A$, then $\forall\varepsilon>0$, $\exists a\in A_{sa}$ such that $c\leq a\leq \|c\|+\varepsilon$.

First assume $p\leq a$ for some $a$ in $A_{sa}$.
Then clearly $p\in \her_{A^{**}}(A)$.
Thus $\forall \varepsilon > 0$, $\exists a'\in A_{sa}$ such that $p\leq a'\leq 1+\var$.
Therefore $p\leq pa' p$, and hence $\alpha(p)\leq 1+\var$.
Since $\var$ is arbitrary, $\alpha(p)\leq 1$.

Now assume $\alpha(p)=1$.
We will prove $p\in\her_{A^{**}}(A)$.
Let $H$ be the Hilbert space of the universal representation of $A$, so that $A^{**}$ is the von Neumann algebra generated by $A$ in $B(H)$.
Represent elements of $A^{**}$ as $2\times 2$ operator matrices relative to $H=pH\oplus (1-p)H$.
Choose $\var > 0$ and $a$ in $A_+$ such that $\|a\| < 1+\var$ and $p\leq pap$.
Let $a=\pmatrix x&y\\ y^*& z\endpmatrix$.
Since
$$
\pmatrix x&y\\ y^*&z\endpmatrix \leq \pmatrix 1+\var&0\\ 0&1+\var\endpmatrix ,\ \pmatrix 1+\var-x&-y\\ -y^*&1+\var-z\endpmatrix \geq 0.
$$
Therefore $\|y\|\leq \|1+\var-x\|^{1\over 2} \|1+\var-z\|^{1\over 2}\leq e^{1\over 2}(1+\var)^{1\over 2}$, since $x\geq 1$ and $z\geq 0$.
Since 
$$
\pmatrix x&0\\ 0&z\endpmatrix\leq \pmatrix x+\|y\|&y\\ y^*&z+\|y\|\endpmatrix,
$$
$p\leq a+\|y\|$.
Let $(e_i)_{i\in D}$ be an approximate identity of $A$.
Then $\limsup \|(1-e_i) p(1-e_i)\|\leq \limsup \| (1-e_i) a(1-e_i)\|+\|y\|\leq \var^{1\over 2} (1+\var)^{1\over 2}$.
Since $\var$ is arbitrary, $\lim \|(1-e_i) p(1-e_i)\|=0$.
This implies $p\in\her_{A^{**}} (A)$.

We review some known facts about pairs of projections.
A complete classification of these, up to unitary equivalence, was given by Dixmier [12].
See also [14], [17], and [19].
If $p$ and $q$ are projections in $B(H)$ with ranges $M$ and $N$, let $H_{11}=M\cap N$, $H_{10}=M\cap N^\perp$, $H_{01}=M^\perp\cap N$, $H_{00}=M^\perp\cap N^\perp$, and $H_0=(H_{11}\oplus H_{10}\oplus H_{01}\oplus H_{00})^\perp$.
A simple example of a pair of projections occurs when $H$ is two dimensional and $p=\pmatrix 1&0\\ 0&0\endpmatrix$, $q=\pmatrix \cos^2\theta&\cos\theta\sin\theta\\ \cos\theta\sin\theta&\sin^2\theta\endpmatrix$ for some $\theta$ in $(0,{\pi\over 2})$.
In the most general example, $(H_0,p|_{H_0},q|_{H_0})$ is a direct integral of such two dimensional examples, for various values of $\theta$.
$\|p-q\|$ can be computed as follows:\ If $H_{01}$ or $H_{10}$ is non-trivial, $\|p-q\|=1$.
Otherwise $\|p-q\|=\sin\theta$, where $\theta$ is the essential supremum of the angles occurring in the decomposition of $H_0$.
For later use, we make a couple of other points:

1.\ The usual concept of the angle between two projections (or subspaces) is the essential infimum of the angles occurring in the decomposition of $H_0$.

2.\ Define $d_a(p,q)=\sin^{-1} (\|p-q\|)$.
Then $d_a$ is a metric on the set of projections, equivalent to the metric induced by the norm ([9, Corollary 4]).
\enddemo

\proclaim{Theorem 2.2}If $p$ is a projection in $A^{**}$, then dist$(p,RC)=[1-\alpha(p)^{-1}]^{1\over 2}$.
\endproclaim

\demo{Proof}1.\ dist$(p,RC)\geq [1-\alpha(p)^{-1}]^{1\over 2}$:\ For this we may assume dist$(p,RC)<1$.
Let $q$ be in $RC$ such that $\|p-q\| < 1$.
Then $pqp > (\cos^2\theta)p$, where $\theta$ is as above, so that $\|p-q\|=\sin\theta$.
Since $q$ is relatively compact, there is $a$ in $A_{sa}$ such that $q\leq a\leq 1$.
Thus $pap\geq (\cos^2\theta) p$, and
$$
\gathered
\cos^{-2} \theta\geq\alpha(p).\\
\text{Therefore }\cos^2\theta\leq\alpha (p)^{-1},\\
\sin^2\theta\geq 1-\alpha(p)^{-1},\\
\|p-q\|=\sin\theta\geq [1-\alpha (p)^{-1}]^{1\over 2}.
\endgathered
$$
Since $q$ can be chosen so that $\|p-q\|$ approximates dist$(p,RC)$, we conclude that dist$(p,RC)\geq [1-\alpha (p)^{-1}]^{1\over 2}$.

2.\ dist$(p,RC)\leq [1-\alpha(p)^{-1}]^{1\over 2}$:\newline
For this we may assume $\alpha(p) < \infty$.
Let $a$ be in $A_{sa}$ such that $p\leq pap$, let $\var > 0$, and let $q=E_{[\var,\infty)}(a)$($q$ is a spectral projection of $a$).
Then $q$ is compact.
Since $a\leq \|a\|q+\var(1-q)$, $p\leq \|a\|pqp+\var p$.
Therefore $pqp\geq {1-\var\over \|a\|} p$.
Let $r$ be the range projection of $qp$.
Then $r\leq q$ and hence $r\in RC$.
Since $rp=qp$, $prp=(rp)*(rp)=(qp)*(qp)=pqp$.
Refer to the notation introduced above for the pair $(p,r)$.
If $\var < 1$, the initial projection of $rp$ is $p$, and hence $H_{10}=0$.
Since $r$ is the range projection of $rp$, $H_{01}=0$.
Therefore $\|p-r\|=\sin\theta$, where $\cos^2\theta\geq {1-\var\over \|a\|}$.
Thus dist$(p,RC)\leq [1-{1-\var\over \|a\|}]^{1\over 2}$.
We can choose $\var$ and $a$ so that ${1-\var\over \|a\|}$ approximates $\alpha(p)^{-1}$.
\enddemo

\proclaim{Corollary 2.3}$\alpha(p)=1$ if and only if $p$ is in the norm closure of $RC$.
\endproclaim

\remark{Remark}Clearly if $\alpha(p)<\infty$ and $p'$ is sufficiently close to $p$, then $\alpha(p')<\infty$.
If one wants the best estimates (i.e., \underbar{how} close must $p'$ be to $p$ and what is the best estimate for $\alpha(p')$?), one should use the metric $d_a$.
Thus $d_a(p,RC)=\cos^{-1}(\alpha(p)^{-{1\over 2}})$; and if $d_a(p',p)+d_a(p,RC) < {\pi\over 2}$, then $\alpha(p')<\infty$.
It is easy to construct examples (Section 3) where $d_a(p',p)+d_a(p,q)={\pi\over 2}$, $q\in RC$, and $\alpha(p')=\infty$. 
\endremark

\proclaim{Proposition 2.4}Let $h$ be a strongly usc element of $A_+^{**}$ such that the spectrum of $h$ omits $(0,\var)$ for some $\var>0$.
Then $E_{(0,\infty)}(h)$ is compact.
\endproclaim

\demo{Proof}Proposition 2.44(b) of [7] asserts that $E_{(0,\infty)}(h)$ is closed under the hypothesis that $h$ is weakly usc.
The proof of the present result is almost identical.
Alternatively, the present result can be deduced from the earlier one by adjoining an identity to $A$.
\enddemo

\proclaim{Lemma 2.5}Assume $p$ is a projection, $0\leq a\leq 1$, and $pap\geq \var p$ for some $\var>0$.
Then $pa^{1\over 2} (pap)^{-1} a^{1\over 2} p\geq pap$, where the inverse is taken in $pA^{**} p$.
\endproclaim

\remark{Remark}Of course this is an operator-theoretic lemma that has nothing to do with $A$.
\endremark

\demo{Proof}Again we represent elements of $A^{**}$ as $2\times 2$ operator matrices relative to $H=ph\oplus (1-p)H$.
Write $a^{1\over 2}=\pmatrix x&y\\ y^*&z\endpmatrix$, so that $a=\pmatrix x^2+yy^*& *\\ *&*\endpmatrix$\newline
Since $a^{1\over 2}\geq a$, $x\geq x^2+yy^*$.
Therefore\newline
$x(x^2+yy^*)^{-1} x\geq x(x)^{-1} x=x\geq x^2+yy^*$.
This is the desired inequality.
\enddemo

\proclaim{Theorem 2.6}Let $p$ be a projection in $A^{**}$.
\item{(a)}If $p$ is open, then dist$(p,ORC)=$ dist$(p,RC)$.
\item{(b)}If $p$ is closed, then dist$(p,CRC)=$ dist$(p,RC)$.
Moreover, in this case, if $\exists a\in A_{sa}$ such that $p\leq pap$ and $\|a\|=\alpha(p)$, then $\exists q\in CRC$ such that $\|p-q\|=$ dist$(p,CRC)$.
\endproclaim

\demo{Proof}The proofs of the two cases are similar.
We start with $a$ in $A_{sa}$ such that $0\leq a\leq 1$ and $pap\geq\var p$, $\var>0$, and let $q$ be the range projection of $a^{1\over 2}p$.
$\var$ should approximate (or, for the last sentence of (b), be equal to) $\alpha(p)^{-1}$.
In case (a), we also need the range projection of $a$ to be in $RC$.
This is accomplished by replacing the original $a$\newline
with $f_\delta(a)$, where $f_\delta(t)=\cases 0,\ 0\leq t\leq\delta\\ t,\ 2\delta\leq t\leq 1\endcases$.\newline
(This causes the original $\var$ to be replaced by $\var-2\delta$.)

The partial isometry in the polar decomposition of $a^{1\over 2}p$ is $u=a^{1\over 2}p(pap)^{-{1\over 2}}$.
Thus $q=uu^*=a^{1\over 2}(pap)^{-1} a^{1\over 2}$.
Lemma 2.5 implies $pqp\geq pap\geq\var p$.
Also $qpq=a^{1\over 2}(pap)^{-1} (pa^{1\over 2} p)^2 (pap)^{-1} a^{1\over 2}$.
Since $pap$ and $pa^{1\over 2}p$ are invertible elements of $pA^{**}p$, this implies $qpq\geq\delta_1(a^{1\over 2} pa^{1\over 2})$, for some $\delta_1 > 0$, and hence $qpq\geq\delta_2 q$ for some $\delta_2 > 0$.
Thus the range projection of $qp$ is $q$.
Now the discussion preceding and the proof of Theorem 2.2 imply that $\|p-q\|\leq (1-\var)^{1\over 2}$.

To complete the proof, we need only show that $q$ is in $ORC$ or $CRC$ in the two cases.
$q$ is the range projection of $(a^{1\over 2}p)(a^{1\over 2}p)*=a^{1\over 2} pa^{1\over 2}$.
In case (a), $a^{1\over 2} pa^{1\over 2}$ is strongly lsc, and in case (b), $a^{1\over 2} pa^{1\over 2}$ is strongly usc.
(This follows, for example, from 2.44(a) of [7].)
In case (a), it follows from 2.44(a) of [7] that $q$ is open.
Since $q$ is smaller than the range projection of $a$, which is in $RC$, $q$ is in $ORC$.
In case (b), Proposition 2.4 implies that $q$ is compact.
We need to know that $\sigma(a^{1\over 2} pa^{1\over 2})$ omits $(0,\var)$, and this follows from $\sigma(a^{1\over 2} pa^{1\over 2})\cup \{0\}=\sigma (pap)\cup \{0\}$.
\enddemo

\proclaim{Corollary 2.7}If $p$ is an open projection in $A^{**}$, then $\alpha(p)=1$ if and only if $p$ is in the norm closure of $ORC$.
\endproclaim

\remark{Remark}We also recover (in different language) a result of Akemann [4]:\ If $p$ is closed, then $\alpha(p)=1$ if and only if $p$ is compact.
(see Section 1.)
\endremark

We now consider other interpretations of $\alpha(p)$.
Some of these can be considered as methods of computing $\alpha(p)$.

\proclaim{Proposition 2.8}Let $p$ be a non-zero projection in $A^{**}$, $(e_i)_{i\in D}$ an approximate identity of $A$, and $\var_i$ the least point in $\sigma(pe_i p)$, where the spectrum is computed in $pA^{**}p$.
Then $\alpha(p)^{-1}=\lim\var_i$.
\endproclaim

\demo{Proof}$\var_i\leq\alpha (p)^{-1}$ so that $\limsup \var_i\leq\alpha(p)^{-1}$.
(We do not need to assume that $(e_i)$ is increasing, though we do assume $0\leq e_i\leq 1$.)

Assume $0\leq a\leq 1$ and $pap\geq\var p$.
For any $\delta>0$, there is $i_0$ such that $\|a-e_i ae_i\| < \delta$ for $i\geq i_0$.
Thus $\var p\leq pap\leq p(e_i ae_i+\delta)p\leq pe_i^2 p+\delta p\leq pe_i p+\delta p$.
Therefore $\var-\delta\leq \var_i$ for $i\geq i_0$, and $\liminf\var_i\geq\var$.
Since $\var$ can be chosen to approximate $\alpha(p)^{-1}$, $\liminf \var_i\geq \alpha (p)^{-1}$.
\enddemo

\remark{Remark}It was pointed out in [8] (Remark 1 after Theorem 4) that if $e$ is a strictly positive element of $A$, then $\alpha(p) < \infty$ if and only if $pep\geq\var p$ for some $\var>0$.
\endremark

\proclaim{Theorem 2.9}Let $p$ be a non-zero projection in $A^{**}$, and let $\overline S(p)$ be the weak$^*$ closure of $F(p)\cap S(A)$.
Then $\alpha(p)^{-1}=\inf \{\|\varphi\|\colon \varphi\in\overline S(p)\}$.
\endproclaim

\remark{Remarks}1.\ The infimum is actually a minimum.

2.\ This result is most natural when $p$ is closed, but it is valid generally.

3.\ If $p=1$, there is a well known dichotomy:\ If $A$ is unital, $\overline S(1)=S(A)$; and if $A$ is non-unital, $\overline S(1)=Q(A)$.
In our language, $\alpha(1)=1$ or $\infty$ according as $A$ is unital or not.
\endremark

\demo{Proof}Assume $a\in A_{sa}$ and $pap\geq p$.
Then $\varphi(a)\geq 1$, $\forall \varphi\in F(p)\cap S(A)$.
Therefore $\varphi(a)\geq 1$, $\forall\varphi\in\overline S(p)$.
Thus $\|a\|\geq \|\varphi\|^{-1}$, $\forall \varphi\in\overline S(p)$.
This implies $\alpha(p)^{-1}\leq\inf \{\|\varphi\|\colon \varphi\in \overline S(p)\}$.

To prove the reverse inequality, we may assume $\inf\{\|\varphi\|\colon\varphi\in\overline S(p)\} > 0$.
Choose $\var$ such that $0<\var <\inf \{\|\varphi\|\colon \varphi\in \overline S(p)\}$, and let $K=\{f\in A^*\colon f=f^*$ and $\|f\|\leq\var\}$.
Then $K$ and $\overline S(p)$ are disjoint compact convex sets.
By the separation theorem, we can find $a$ in $A_{sa}$ such that $\sup\{ f(a)\colon f\in K\} < \inf \{\varphi (a)\colon \varphi\in \overline S(p)\}$.
Since the supremum is $\var \| a\|$, we can normalize $a$ so that $\|a\|=1$, and then we find $pap\geq \var p$.
This implies $\alpha(p)^{-1}\geq\var$ and hence $\alpha(p)^{-1}\geq \inf \{\|\varphi\|\colon \varphi\in \overline S(p)\}$.
\enddemo

\proclaim{Corollary 2.10}$\alpha(p) < \infty$ if and only if 0 is not in the weak$^*$ closure of $F(p)\cap S(A)$.
\endproclaim

If $V$ is a partially ordered real normed linear space and $e\in V_+$, $e$ is an \underbar{order unit} of $V$ if $\forall x\in V$, $\exists t\in\bR_+$ such that $x\leq te$.
We will call $e$ a \underbar{$t$-order unit} if $\|e\|=1$ and $x\leq t\|x\| e$, $\forall x\in V$.
If $V$ is a Banach space and the positive cone is closed, then every order unit of norm 1 is a $t$--order unit for $t$ sufficiently large.
The proof of this (presumably known) result is similar to an argument given in the next theorem.
If $p$ is a projection in $A^{**}$, then $pA_{sa} p$ is a partially ordered real normed linear space if regarded as a subspace of $pA_{sa}^{**}p$.
If $p$ is closed, then a result of [6] implies that $pA_{sa}p$ is a Banach space and its norm is the quotient norm from the natural map $A_{sa}\to p A_{sa} p$.

\proclaim{Theorem 2.11}Let $p$ be a projection in $A^{**}$.\newline
(a)\ $\alpha(p) < \infty$ if and only if $pA_{sa}p$ has an order unit.\newline
(b)\ If $p$ is closed, then $\alpha(p)=\inf\{t\colon pA_{sa}p$ has a $t$--order unit$\}$.
Also $pA_{sa}p$ has an $\alpha(p)$-order unit if and only if there is $a$ in $A_{sa}$ such that $pap\geq p$ and $\|a\|=\alpha(p)$\newline
(c)\ For general $p,\alpha(p)$ is the infimum of $t$ such that there is an order unit $e$ satisfying:
\item{(i)}$e=pa_1 p$ where $\|a_1\|\leq 1$.
\item{(ii)}$pap\leq t\|a\|e,\ \forall a\in A_{sa}$.
\endproclaim

\demo{Proof}(a)\ If $pap\geq p$, then clearly $pap$ is an order unit for $pA_{sa}p$.
Conversely, if $e$ is an order unit, let\newline
$C=\{a\in A_{sa}\colon -e\leq pap\leq e\}$.
Then $C$ is closed, convex, and symmetric, and $A_{sa}=\bigcup\limits^\infty_1 n C$.
A standard argument based on the Baire category theorem shows that $nC$ contains the unit ball of $A_{sa}$, for some $n$.
If $(e_i)$ is an approximate identity of $A$, then $pe_i p\leq ne$, $\forall i$.
Taking the strong limit, we see that $p\leq ne$.
Therefore $\alpha(p) < \infty$.

(c)\ If $e$ and $t$ satisfy (i) and (ii), then part of the argument just given shows that $p\leq te=p(ta_1)p$.
Thus $\alpha(p)\leq t$.
Therefore $\alpha(p)$ is at most the infimum specified.
On the other hand, if $pap\geq p$, then $e$ and $t$ satisfy (i) and (ii), where $t=\|a\|$ and $e=t^{-1} pap$.
This implies the opposite inequality.

(b)\ If $p$ is closed, the infima in (b) and (c) are the same, since the norm of $pA_{sa}p$ is the quotient norm under the map $a\mapsto pap$ ([6]).
The second sentence of (b) is deduced from 3.3 or 3.4 of [7]:\ $e=pap$ where $\|a\|=\|e\|$.
\enddemo

\proclaim{Lemma 2.12}Assume $p$ is a closed projection in $A^{**}$, $a\in A_+$, and $pap\geq\var p$ for some $\var>0$.
Then $pa^{1\over 2} Aa^{1\over 2}p=pAp$.
\endproclaim

\demo{Proof}Since $pap\geq\var p$, $\exists s\in A^{**}$ such that $p=sa^{1\over 2}p$.
Then $pA^{**} p=pa^{1\over 2} (s^* A^{**} s) a^{1\over 2} p\subset pa^{1\over 2} A^{**} a^{1\over 2} p$.
Now let $y=pa^{1\over 2}$, and define $\varphi\colon A\to pAp$ by $\varphi(b)=yby^*$.
The result cited from [6] shows that $(pAp)^{**}$ can be identified with $pA^{**}p$ in such a way that $\varphi^{**}$ becomes the map of $A^{**}$ to $pA^{**}p$ given by $b\mapsto yby^*$.
In general, if the second adjoint of a map between Banach spaces is surjective, then the original  map is surjective.
\enddemo

\proclaim{Theorem 2.13}Let $p$ be a closed projection in $A^{**}$.\newline
(a)\ $\alpha(p)<\infty$ if and only if there are a compact projection $q$ and a complete order isomorphism $\theta\colon pAp\to qAq$.\newline
(b)\ $\alpha(p)=\inf\{\|\theta\| \|\theta^{-1}\|\colon q$ and $\theta$ as above$\}$.\newline
(c)\ There are $\theta$ and $q$ as above such that $\|\theta\| \|\theta^{-1}\|=\alpha(p)$ if and only if there is $a$ in $A_{sa}$ such that $pap\geq p$ and $\|a\|=\alpha(p)$.
\endproclaim

\demo{Proof}First assume $\theta$ and $q$ are as in (a).
Since $q$ is compact, $q\in qAq$.
Let $e=\theta^{-1}(q)$.
If $b\in A_{sa}$, then $\theta(pbp)\leq \|\theta\| \|b\|q$.
Therefore $pbp\leq \|\theta\| \|b\|e$.
As in the proof of theorem 2.11, we deduce that $p\leq \|\theta\|e$.
By 3.4 of [7], we can write $\|\theta\|e=pap$ for $a$ in $A_{sa}$ such that $\|a\|=\|\theta\| \|e\| \leq \|\theta\| \|\theta^{-1}\|$.

Next assume $\alpha(p) < \infty$, $a\in A$, $0\leq a\leq 1$, and $pap\geq\var p$.
Here $\var$ approximates $\alpha(p)^{-1}$, and for (c), $\var=\alpha(p)^{-1}$.
Let $q$ be the range projection of $a^{1\over 2}p$.
As in the proof of Theorem 2.6, we deduce that $q$ is compact and $q=a^{1\over 2}(pap)^{-1} a^{1\over 2}$.
Then \newline
$qAq=a^{1\over 2} (pap)^{-1} (a^{1\over 2} Aa^{1\over 2}) (pap)^{-1} a^{1\over 2}= a^{1\over 2} (pap)^{-1} A(pap)^{-1} a^{1\over 2}$,\newline
where the second equality uses Lemma 2.12.
Let $x=a^{1\over 2} (pap)^{-1}$ and $y=pa^{1\over 2}$.
Then $x\in qA^{**}p,\ y\in pA^{**} q,\ xy=q$, and $yx=p$.
If we define $\theta$ and $\varphi$ by $\theta(b)=xbx^*$, and $\varphi(b)=yby^*$, then the above equation shows that $\theta$ maps $pAp$ into $qAq$ and it is obvious that $\varphi$ maps $qAq$ into $pAp$.
It is now obvious that $\theta$ and $\varphi$ are inverses of one another, and clearly both are completely positive.
Now\newline
$\|\theta\| \leq \|x\|^2=\|x^* x\|=\|(pap)^{-1}\|\leq \var^{-1}$, and $ \|\theta^{-1}\|\leq \|y\|^2\leq 1$.
\enddemo

The above arguments prove all three parts of the theorem.

\remark{Remark}The first part of the proof used only the hypothesis that $\theta$ is an order isomorphism, not a \underbar{complete} order isomorphism.
Therefore, the word ``complete'' could be omitted from the statement of the theorem.
\endremark

\subheading{\S 3.\ Some Examples and Discussion}

If $1\leq s\leq t\leq \infty$, $A$ is a $C^*$-algebra, $p$ is a projection in $A^{**}$, and $(\alpha(p),\alpha(\overline p))=(s,t)$, we will say that $p$ and $A$ achieve $(s,t)$, or, more briefly, that $p$ achieves $(s,t)$.
The basic object of this section is to show that every such pair can be achieved, but we want a little more.
We want to consider various properties of projections and find which pairs can be achieved by projections satisfying one or more of these properties.
The properties we will consider are open, closed, central, and regular, except that all discussion of regularity will be postponed to the next section (this does not cause much inefficiency).
Of course there are many other properties which could be considered, and perhaps some of these would lead to deeper results.
The gist of what we will show is that all pairs can be achieved with open projections, but the other properties are compatible only with very special pairs.
If $p$ is closed, obviously we must have $s=t$.
The restrictions required for the other properties are not much deeper, but we will dignify them with numbers.

\noindent
{\bf 3.1}.
If $p$ is clopen (both open and closed), then either $\alpha(p)=\alpha(\overline p)=1$ or $\alpha(p)=\alpha(\overline p)=\infty$.

\demo{Proof}Of course $p$ is clopen if and only if $p\in M(A)$, the multiplier algebra of $A$.
If $pap\geq p$ for $a$ in $A$, then $pap$ is in $A$ also.
From this we easily conclude that $p$ is in $A$ (look at the images in $M(A)/A$).
\enddemo

\noindent
{\bf 3.2}.
If $p$ is a central projection in $A^{**}$, then either $\alpha(p)=\alpha(\overline p)=1$ or $\alpha(p)=\alpha(\overline p)=\infty$.

\demo{Proof}If $\alpha(p) <\infty$, then there is $a$ in $A_+$ such that $pap\geq p$.
Since $pa=ap$, this clearly implies $a\geq p$.
Let $f(x)=\cases x,\quad 0\leq x\leq 1\\ 1,\quad x\geq 1\endcases$.
Since $a\geq p$ and $ap=pa$, $1\geq f(a)\geq p$.
Therefore $\alpha(\overline p)=1$.

In the examples $\cK$ denotes the set of compact operators on a separable infinite dimensional Hilbert space $H$, $ \{e_1,e_2,\ldots\}$ is an orthonormal basis of $H$, and $ v\times w$ denotes the rank 1 operator $x\mapsto (x,w)v$.
In many cases we will take $A=c\otimes\cK$.
Then $A$ can be regarded as the set of $\{x_n\colon 1\leq n\leq\infty\}$ such that $x_n\in\cK$ and $x_n\to x_\infty$ in norm, and $A^{**}$ is the set of $\{h_n\colon 1 \leq n \leq \infty\}$ such that $h_n\in B(H)$ and $\{\|h_n\|\}$ is bounded.
If $p=\{p_n\}$ is a projection in $(c\otimes\cK)^{**}$, then $p$ is open if and only if $p_\infty\leq h$ for every weak cluster point $h$ of the sequence $(p_n)$, and $p$ is closed if and only if $p\geq h$ for every such $h$.
This follows, for example, from the criterion for weak semicontinuity given in 5.14 and 5.15 of [7].
\enddemo

\noindent
{\bf 3.3}.\ (1,1).

It is trivial to achieve this pair with a clopen central projection.
Just let $A$ be any unital $C^*$-algebra and $p=1$.

\noindent
{\bf 3.4}.
$(\infty,\infty)$.
It is trivial to achieve this pair with a clopen central projection.
Just let $A$ be any non-unital $C^*$-algebra and $p=1$.

\noindent
{\bf 3.5}.
$(s,s),\ 1< s<\infty$ (cf.~Remark 2, p.~276, of [8]).

For this we need two examples, one open and one closed.
Let $\theta$ be in $(0,{\pi\over 2})$ such that $\sec^2\theta=s$.
Let $A=c\otimes\cK$ and $v_n=\cos\theta e_1+\sin\theta e_{n+1}$.
Define $p$ and $q$ by $p_n=q_n=v_n\times v_n, n<\infty$, $p_\infty=e_1\times e_1$, and $q_\infty=0$.
Then $q$ is open, $p=\overline q$, and we claim $\alpha(p)=\alpha(q)=s$.
(Thus $p$ and $q$ both achieve the pair $(s,s)$).
Define $a$ in $A$ by $a_n=a_\infty=s(e_1\times e_1)$.
Then $pap\geq p$ and $qaq\geq q$ (actually $qaq=q$).
Thus $\alpha(q), \alpha(p)\leq s$.
If $\varphi_n$ is the pure state of $A$ given by $\varphi_n(a)=(a_n v_n,v_n)$ then $\varphi_n\in F(q)\cap S(A)\subset F(p)\cap S(A)$ and $\varphi_n$ converges weak$^*$ to a functional of norm ${1\over s}$.
Thus $\alpha(q)$, $\alpha(p)\geq s$.

We now justify the remark after Corollary 2.3.
Choose $\theta'$ such that $\theta<\theta'\leq {\pi\over 2}$, and let $w_n=\cos\theta' e_1+\sin\theta' e_{n+1}$.
Define a closed projection $p'$ by $p'_\infty=e_1\times e_1$, $p'_n=w_n\times w_n$.
Then $\alpha(p')=\sec^2\theta'$, $d_a(p',RC)=\theta'$, $d_a(p',p)=\theta'-\theta$, and $d_a(p,RC)=\theta$.
If $\theta'={\pi\over 2}$, $\alpha(p')=\infty$.
We could equally well consider an open projection $q',\ (q')_\infty=0$, $q'_n=p'_n$, and compare $q'$ to $q$.

\noindent
{\bf 3.6}.
$(s,\infty)$, $1<s<\infty$ (cf.~Remark 4, p.~276, of [8]).

Let $A,\theta$, and $v_n$ be as in the previous example.
Let $(m_n)$ be a sequence which includes each positive integer infinitely often.
Define an open projection $p$ in $A^{**}$ by $p_\infty=0$, $p_n=v_{m_n}\times v_{m_n}$, $n<\infty$.
Then by essentially the same argument as above, $\alpha(p)=s$.
Since $\{v_n\}$ is total in $H$, $(\overline p)_\infty=1$.
Then it is obvious that $\alpha(\overline p)=\infty$.

\noindent
{\bf 3.7}.\ $(1,\infty)$.

Akemann's [4, Example IV.5] gives an open projection that achieves this pair, but we will give another example, somewhat similar in spirit, where $A=c\otimes\cK$.

\proclaim{Lemma}If $x > 1 > y > 0$ and $u=(u_1,u_2)$ where $|u_1|^2={x(1-y)\over x-y}$ and $|u_2|^2={y(x-1)\over x-y}$, then $u\times u\leq \pmatrix x&0\\ 0&y\endpmatrix$.
(Here the Hilbert space is 2-dimensional.)
\endproclaim

Now let $K$ be any one-to-one element of $\cK_+$ such that $\|K\| > 1$.
If $\cV=\{ u\in H\colon \|u\|=1$ and $u\times u\leq K\}$, then the lemma implies that $\cV$ is a total subset of $H$.
Let $(u_n)$ be a sequence which is dense in $\cV$.
Define an open projection $p$ in $A^{**}$ by $p_\infty=0$, $p_n=u_n\times u_n$.
Define $a$ in $A$ by $a_n=a_\infty=K$.
Since $p\leq a$, $\alpha(p)=1$.
Since $(\overline p)_\infty=1$, $\alpha(\overline p)=\infty$.

\noindent
{\bf 3.8}.
$(1,t), 1< t<\infty$.

Let $A_0=c\otimes\cK$ and $p_0$ the projection in $A_0^{**}$ called $p$ in 3.7.
Let $A_1=A_0\otimes M_2$.
For this example, $A$ will be an extension of $A_1$ by $\bC$.
According to Busby [10], such an extension is determined by an element $e'$ of $M(A_1)$ which maps onto a projection in $M(A_1)/A_1$.
We will take $e'$ to actually be a projection; namely,
$e'=\pmatrix t^{-1}&[t^{-1}(1-t^{-1})]^{1\over 2}\\ [t^{-1}(1-t^{-1})]^{1\over 2}& 1-t^{-1}\endpmatrix$.
Let $e$ be the corresponding element of $A$.
Thus $e^2=e=e^*$ and $ex=e'x$, $xe=xe'$ for $x$ in $A_1$.
Then $A^{**}\simeq A_1^{**}\oplus \bC\simeq (A_0^{**}\otimes M_2)\oplus\bC$.

Now let $p=\pmatrix p_0&0\\ 0&0\endpmatrix \oplus 0$.
Clearly $p$ is open.
If $a$ has the same meaning as in 3.7 (so that $a\in A_0$), then $p\leq\pmatrix a&0\\ 0&0\endpmatrix\oplus 0$, an element of $A$.
Thus $\alpha(p)=1$.
We claim that $\overline p=\pmatrix \overline p_0&0\\ 0&0\endpmatrix \oplus 1$.
In fact, clearly $\overline p=\pmatrix \overline p_0&0\\ 0&0\endpmatrix \oplus r$, where $r$ is 0 or 1.
It is actually not important which is true.
To show that $r=1$, we need only show $\not\exists x\in A_1$ such that $(e-x)p=0$.
This is equivalent to showing $\not\exists x'\in A_0$ such that $(1-x') p_0=0$.
This last follows from the fact that $\overline p_0$ is not compact in $A_0^{**}$.

Obviously $\overline p\leq \overline p(te)\overline p$.
Therefore $\alpha(\overline p)\leq t$.
Since $\alpha(\overline p_0=\infty$, there is a sequence $(\varphi'_n)$ in $F(\overline p_0)\cap S(A_0)$ such that $\varphi'_n\to 0$ in the weak$^*$ topology of $A^*_0$.
Let $(\varphi_n)$ be the corresponding sequence in $F(\overline p)\cap S(A)$.
Note that $A^*\simeq A^*_1\oplus\bC$ and $\varphi_n\in A_1^*\oplus 0$.
Since $\varphi'_n\to 0$, every weak$^*$ cluster point of $(\varphi_n)$ has $A_1^*$-component 0.
Since $\varphi_n(e)=t^{-1}$, $\forall n$, we conclude that $\varphi_n\to 0\oplus t^{-1}$ in the weak$^*$ topology of $A^*$.
Therefore $\alpha(\overline p)\geq t$.

\noindent
{\bf 3.9}.\ $(s,t)$, $1< s<t<\infty$.

\remark{Remark}If one is only interested in which of (1), (2), (3) (notation of Section 1) are satisfied by $p$ and $\overline p$, then it is not necessary to consider this example, since 3.5 would suffice.
\endremark

For this example $A$ is the same as in 3.8.
In particular $A_0,A_1$, $e'$ and $e$ are the same.
Let $p_0$ be the projection in $A_0^{**}$ called $p$ in 3.6, with the $s$ of 3.6 replaced by $s'$, where $s'$ is a number in $(1,\infty)$ to be determined later.
As in 3.8, we let $p=\pmatrix p_0&0\\ 0&0\endpmatrix \oplus 0$, an open projection in $A^{**}$, and $\overline p=\pmatrix \overline p_0&0\\ 0&0\endpmatrix \oplus 1$.

As in 3.8, we prove that $\alpha(\overline p)=t$.
It remains to calculate $\alpha(p)$.
Since $\alpha (p_0)=s'$, there is a sequence $(\varphi'_n)$ in $F(p_0)\cap S(A_0)$ such that $\varphi'_n\to\varphi'$, where $\|\varphi'\|=(s')^{-1}$, in the weak$^*$ topology of $A_0^*$.
Let $(\varphi_n)$ be the corresponding sequence in $F(p)\cap S(A)$, and let $\varphi$ be the element of $A_1^*$ corresponding to $\varphi'$.
Since $\varphi_n(e)=t^{-1}$, $\forall n$, we find that $\varphi_n\to \varphi\oplus (t^{-1} - (s')^{-1} t^{-1})$ in the weak$^*$ topology of $A^*$ (cf.~3.8).
Since $\|\varphi\oplus (t^{-1} -(s')^{-1} t^{-1}\|=(s')^{-1} + t^{-1} - (s')^{-1} t^{-1}=(s')^{-1}+(1-(s')^{-1})t^{-1}$, $\alpha(p)^{-1}\leq (s')^{-1}+(1-(s')^{-1})t^{-1}$.
Let $q$ be the projection in $A_0$ given by $q_n=q_\infty=e_1\times e_1$.
Then $p_0 qp_0=(s')^{-1} p_0$.
Then we define an element $a$ of $A$ by
$$
a=e+\pmatrix (1-t^{-1})q& -[t^{-1} (1-t^{-1})]^{1\over 2} q\\ -[t^{-1}(1-t^{-1})]^{1\over 2}q & -(1-t^{-1})q\endpmatrix,
$$
where the matrix is in $A_1$.
Thus, relative to $A^{**}\simeq (A_0^{**}\otimes M_2)\oplus \bC$, we have $a=\pmatrix q+t^{-1}(1-q) & [t^{-1} (1-t^{-1})^{1\over 2}](1-q)\\ [t^{-1} (1-t^{-1})]^{1\over 2} (1-q)& (1-t^{-1})(1-q)\endpmatrix \oplus 1$.
Clearly, $\|a\|=1$, and $pap=[(s')^{-1}+(1-(s')^{-1})t^{-1}]p$.
Now, if we choose $s'$ such that $(s')^{-1}+(1-(s')^{-1})t^{-1})=s^{-1}$, we have that $\alpha(p)=s$.

\subheading{\S 4.\ Regularity, Some Variants, and Relations with $\alpha(p)$}

Before proceeding, the author has to make a personal statement:\ In 1985 I was told that someone had done some work on variants of regularity.
Specifically, I was told this mathematician's definition of $k$-regularity (given below); and I think I was told there was a special result on 2-regularity, but I was not told what this result was (it is likely similar to my 4.16, 4.17).
Unfortunately, I was not interested enough then to ask this mathematician's name, and now (1990) the person who told me has forgotten the name.
I made a strong effort to locate a name or paper without success.
Except as noted above all of my work is independent, in particular all of my proofs are independent, but surely some of my results were obtained first by the inventor of $k$-regularity.
Except for one comment in Example 4.15(b), I make no further reference to this unpleasant situation.

For $p$ a projection in $A^{**}$, we have already defined $F(p)$, the norm closed face of $Q(A)$ supported by $p$.
There are many other convex subsets of $A^*$ that can be defined in terms of $p$.
Among these:\ $L(p)=\{f\in A^*\colon f(a)=f(ap),\ \forall a\}$, the left ideal generated by $F(p)$.
$\{L_1(p)=\{f\in L(p)\colon \|f\|\leq 1\}$.
$C(p)=\{f\in A^*\colon f\geq 0$ and $f(1-p)=0\}$, the cone generated by $F(p)$.
$V(p)=\{f\in A^*\colon f(a)=f(pap),\ \forall a\}$, the complex vector space generated by $F(p)$.
$RV(p)=\{f\in V(p)\colon f=f^*\}$, the real vector space generated by $F(p)$.
$V_1(p)=\{f\in V(p)\colon \|f\|\leq 1\}$.
$RV_1(p)=\{f\in RV(p)\colon \|f\|\leq 1\}$, the convex hull of $F(p)\cup (-F(p))$.

If $p$ is closed, then all of the above sets are weak$^*$ closed; and if any of these sets is weak$^*$ closed, then $p$ is closed.
All of these facts were either proved by Effros in [13] or are easy consequences of results of [13].
The problem of relating the closure operation to these sets is more complicated.
Effros showed that $L(p)^- =L(\overline p)$, where ``$-$'', when applied to a subset of $A^*$, always means weak$^*$ closure.
We will use the following uninspired abbreviations:

\noindent
$(R_1)\ \ L_1(p)^- = L_1 (\overline p)$\newline
$(R_2)\ \ F(p)^- = F(\overline p)$\newline
$(R_3)\ \ C(p)^- =C(\overline p)$\newline
$(R_4)\ \ V(p)^- =V(\overline p)$\newline
$(R'_4)\ \ RV(p)^- = RV(\overline p)$\newline
$(R_5)\ \ V_1 (p)^- = V_1(\overline p)$\newline
$(R_6)\ \ RV_1 (p)^- =RV_1 (\overline p)$\newline
$(R_7 (K))\ L_1 (p)^- \supset K^{-1} L_1 (\overline p),\ 1< K<\infty$\newline
$(R_8 (K))\ F(p)^- \supset K^{-1} F(\overline p),\ 1< K<\infty$\newline
$(R_9 (K))\ \ RV_1 (p)^- \supset K^{-1} RV_1 (\overline p),\ 1< K<\infty$

$p$ is called \underbar{regular} (Tomita [21]) if $\|a\overline p\|=\|ap\|$, $\forall a\in A$.
Theorem 6.1 of [13] asserts that regularity is equivalent to each of $(R_1),(R_2),(R_3)$.
Unfortunately, the proof tacitly assumed $A$ to be unital in one place, and the theorem is not correct in the non-unital case.
In general, regularity is equivalent to $(R_1)$ and $(R_2)$, and a correct proof of this is contained in [13], but $(R_3)$ may be strictly weaker.

Each of the $(R_i)$'s is a variant of regularity.
There are some deliberate omissions from the list.
Aside from the one the reader has already noticed, we mention in passing a condition intermediate between $(R_8(K))$ and $(R_3)$:\ Every element of $C(\overline p)$ is the weak$^*$ limit of a bounded net from $C(p)$.
The reason for the omissions is not that we are trying to hide anything.
We are simply trying make a reasonable compromise between, on the one hand, presenting the minimum amount of material on regularity indicated by our interest in near relative compactness, and, on the other hand, attempting an exhaustive treatment of the variants of regularity.
(We have not, in fact, done enough research for the latter course.)

The following implications are either obvious or were proved in [13]:
$$
\gathered
\text{regular }\Leftrightarrow (R_1) \Leftrightarrow (R_2)\Rightarrow (R_3), (R_4), (R'_4), (R_6)\\
(R_4)\Leftrightarrow (R'_4)\\
(R_1)\Rightarrow (R_7(K))\\
(R_2)\Rightarrow (R_8(K))\Rightarrow (R_3)\Rightarrow (R'_4)\\
(R_5)\Rightarrow (R_6)\Rightarrow (R_9(K))\Rightarrow (R'_4)\\
(R_8(K))\Rightarrow (R_9(K)).
\endgathered
$$

In particular, regularity implies all except $(R_5)$, and all except $(R_7(K))$ imply $(R_4)$.
$(R_4)$ is therefore interesting, and we will say $p$ is \underbar{0-regular} if $p$ satisfies $(R_4)$.
We say $p$ is \underbar{cone-regular} if it satisfies $(R_3)$, \underbar{$K$-quasi-regular} if it satisfies $(R_7(K))$, and \underbar{quasi-regular} if $K$-quasi-regular for some $K$.
We believe that cone-regularity and quasi-regularity are the most interesting for near relative compactness, but we may have overlooked something.
Finally, $p$ is \underbar{$k$-regular} if $\pmatrix p&&0\\ &\ddots\\ 0&&p\endpmatrix$ is regular in $(A\otimes M_k)^{**}$.

Before finally getting down to business, we need some more notation.
A small amount of semicontinuity theory is used, and we follow the notation of [5].
$\tilde A=A+\bC 1$, where 1 is the identity of $A^{**}$.
For $S\subset A_{sa}^{**}, S^m$ is the set of ($\sigma$-strong) limits of bounded increasing nets from $S$.
$S_m$ is defined similarly with decreasing nets.
``$-$'', when applied to subsets of $A^{**}$, means norm closure.
For example, $(\tilde A_{sa}^m)^-$ is the set of weakly lsc elements of $A^{**}$, and $(A_+^m)^-$ is the set of positive strongly lsc elements.
$M(A)$ is the multiplier algebra of $A$ and $QM(A)$ the space of quasi-multipliers (both are subsets of $A^{**}$).
In all the results of this section $A$ is an arbitrary $C^*$-algebra and $p$ is a projection in $A^{**}$.
The arguments presented below almost include a new proof of the equivalence of regularity, $(R_1)$, and $(R_2)$; but this is not our goal and we officially are assuming this equivalence.

A good way to deal with the $(R_i)$'s is to use the double polar theorem.
For example, $(R_1)$ is equivalent to the statement that $L_1(p)$ and $L_1(\overline p)$ have the same polar in $A$.
It is easy to compute the polars if one remembers that $A^*$ is the predual of the $W^*$-algebra $A^{**}$ and that the polar in $A$ is just the intersection with $A$ of the polar in $A^{**}$.
The polar in $A$ of $L(p)$ is $\{a\in A\colon ap=0\}$.
(Since always $L(p)^- =L(\overline p)$, this tells us that $ap=0\Leftrightarrow a\overline p=0$.
This is often a good ``working definition'' of the closure of a projection.)
The polar in $A$ of $L_1(p)$ is $\{a\in A\colon \|ap\|\leq 1\}$.
The polar in $A_{sa}$ of $F(p)$ is $\{a\in A_{sa}\colon pap\leq p\}$.
(The polar in $A$ is $\{a\in A\colon \text{ Re }pap\leq p\}$.)
The polar in $A_{sa}$ of $C(p)$ is $\{a\in A_{sa}\colon pap\leq 0\}$.
The polar in $A$ of $V(p)$ is $\{a\in A\colon pap=0\}$, and the polar in $A_{sa}$ of $RV(p)$ is $\{a\in A_{sa}\colon pap=0\}$.
The polar in $A$ of $V_1(p)$ is $\{a\in A\colon \|pap\|\leq 1\}$, and the polar in $A_{sa}$ of $RV_1(p)$ is $\{a\in A_{sa}\colon \|pap\|\leq 1\}$.

The following is now obvious:

\proclaim{Proposition 4.1}(a)\ $p$ is $K$-quasi-regular if and only if $\|a\overline p\|\leq K\|ap\|$, $\forall a\in A$.\newline
(b)\ $(R_9(K))\Leftrightarrow \|\overline p a\overline p\|\leq K\| pap\|$, $\forall a\in A_{sa}$.\newline
(c)\ $p$ is 0-regular if and only if $pap=0\Rightarrow \overline p a\overline p=0,\ \forall a\in A$.\newline
(d)\ $p$ is cone-regular if and only if $pap\leq 0\Rightarrow \overline p a\overline p\leq 0,\ \forall a\in A_{sa}$.\newline
(e)\ $(R_5)\Leftrightarrow \|\overline p a\overline p\|=\| pap\|,\ \forall a\in A$.
\endproclaim

Throughout this section, $\sigma(php)$ means the spectrum of $php$ relative to $pA^{**} p$.

\proclaim{Theorem 4.2}(a)\ $p$ is regular if and only if the top points in $\sigma(php)$ and $\sigma(\overline p h\overline p)$ are the same, $\forall h\in (\tilde A_{sa}^m)^-$.

(b)\ $p$ is $K$-quasi-regular if and only if $\|\overline p h\overline p\|\leq K^2\|php\|,\ \forall h\in\overline{A_+^m}$.
\endproclaim

\demo{Proof}(a)\ Since $\|ap\|^2=\|p(a^* a)p\|$, regularity is equivalent to the condition stated for all $h$ in $A_+$.
In general, if $h_i\nearrow h$ in a $W^*$-algebra, the top point in $\sigma(h_i)$ converges to the top point in $\sigma(h)$.
This generalizes the condition to $A_+^m$.
In general, the top point in $\sigma(p(h+\lambda)p)$ is $\lambda+$ top point in $\sigma(php)$.

This generalizes the condition to $\{h\in A_{sa}^{**}\colon \exists\lambda$ with $h+\lambda\in A_+^m\}=\tilde A_{sa}^m$.
In general, if $h_n\to h$ in norm the top point in $\sigma(h_n)$ converges to the top point in $\sigma(h)$.
This generalizes the condition to $(\tilde A_{sa}^m)^-$.

The proof of (b) is similar except that we leave out the step involving translation by $\lambda$.
\enddemo

\proclaim{Corollary 4.3}(a)\ If $p$ is regular, then $\|T\overline p\|=\|Tp\|$ whenever $T^* T\in (\tilde A_{sa}^m)^-$, in particular whenever $T\in QM(A)$.

(b)\ If $p$ is $K$-quasi-regular, then $\|T\overline p\|\leq K\|Tp\|$ whenever $T^*T\in\overline{A_+^m}$, in particular whenever $T$ is a right multiplier of $A$.

(c)\ If $p$ is regular and $h$ is in $QM(A)_{sa}$, in particular if $h$ is in $A_{sa}$ or $M(A)_{sa}$, then both extreme points of $\sigma(php)$ and $\sigma(\overline p h\overline p)$ agree.

(d)\ $(R_9 (K^2))\Rightarrow K$-quasi-regular

(e)\ $(R_6)\Leftrightarrow$ regular

(f)\ $(R_5)\Rightarrow$ regular
\endproclaim

\demo{Proof}For (a) and (b) we just have to quote [7, 4.1].

(c)\ follows from the fact ([5]) that $QM(A)_{sa}=\{h\colon h,-h\in (\tilde A_{sa}^m)^-\}$.

For (d) we use Proposition 4.1 (b) to characterize $(R_9(K^2))$.
By the proof of 4.2, $(R_7(K))$ is equivalent to the restriction of this condition from $A_{sa}$ to $A_+$.

If we let $K\to 1^+$ in (d), we see that $(R_6)\Rightarrow$ regular.
We already knew the converse.
Thus (e) is proved and (f) follows.

It was proved by Pedersen in [18] that if $A$ is unital, then $p$ is regular if and only if $a\geq p\Rightarrow a\geq\overline p,\ \forall a\in A_{sa}$.
His arguments can be generalized:

\proclaim{Theorem 4.4}(a)\ If there is $a$ in $A_{sa}$ such that $a\geq p$, then $p$ is $K$-quasi-regular if and only if $b\geq p\Rightarrow K^2 b\geq\overline p$, $\forall b\in A_{sa}$.

(a$'$)\ If there is $a$ in $A_{sa}$ such that $a\geq p$, then $p$ is regular if and only if $b\geq p\Rightarrow b\geq\overline p$, $\forall b\in A_{sa}$.

(b)\ If general, $p$ is $K$-quasi-regular if and only if $h\geq p\Rightarrow K^2 h\geq\overline p,\ \forall h\in\tilde A_{sa}$ if and only if $h\geq p\Rightarrow K^2 h\geq\overline p$, $\forall h\in M(A)_{sa}$.
If $A$ is $\sigma$-unital, one can add:\newline
$h\geq p\Rightarrow K^2 h\geq\overline p$, $\forall h\in QM(A)_{sa}$.

(b$'$)\ Same as (b) except omit ``$K$'' and ``quasi-''.
\endproclaim

\demo{Proof}(b)\ Assume $p$ is $K$-quasi-regular and $h\geq p$, where $h$ is in $M(A)_{sa}$ or $QM(A)_{sa}$.
Let $\var>0$, and choose $R$ a right multiplier of $A$ such that $R$ is invertible in $A^{**}$ and $R^*R=(h+\var)^{-1}$.
If $A$ is $\sigma$-unital, the existence of $R$ follows from [7, 4.8].
Otherwise, $h\in M(A)$ and we take $R=(h+\var)^{-{1\over 2}}$.
Then $R^{-1} (R^*)^{-1}=h+\var\geq p$,
$$
\gathered
1\geq Rp R^*,\\
\|Rp\|\leq 1,\\
\|R\overline p\|\leq K,\text{ by }4.3(b),\\
R\overline p R^*\leq K^2,\text{ and hence}\\
\overline p\leq K^2 R^{-1} (R^*)^{-1}=K^2 (h+\var).
\endgathered
$$
Since $\var$ is arbitrary, $\overline p\leq K^2 h$.
\enddemo

Now assume $h\geq p\Rightarrow K^2 h\geq\overline p$, $\forall h\in\tilde A_{sa}$.
Assume $x\in A$, $\|xp\|\leq 1$, and $\var > 0$.
Then 
$$
\gathered
\| |x|p\|\leq 1,\\
\| (|x|+\var)p\|\leq 1+\var,\\
(|x|+\var) p(|x|+\var)\leq (1+\var)^2,\text{ and hence}\\
p\leq (1+\var)^2 (|x|+\var)^{-2}.
\endgathered
$$
By hypothesis, $\overline p\leq K^2(1+\var)^2(|x|+\var)^{-2}$.
By reversing some of the above steps, we obtain $\|(|x|+\var)\overline p\|\leq K(1+\var)$.
Then taking limits as $\var\to 0^+$, we obtain $\|x\overline p\|=\| |x|\overline p\|\leq K$.

(a)\ Half of this follows from (b).
Thus assume $b\geq p\Rightarrow K^2 b\geq\overline p$, $\forall b\in A_{sa}$ and $\exists a\in A_{sa}$ such that $a\geq p$.
Clearly then, $\overline p$ satisfies (2) of Section 1, and hence $\overline p$ is compact.
Thus we can choose $a$ in $A_{sa}$ such that $\overline p\leq a\leq 1$.
Now let $(f_i)$ be any approximate identity of $A$, and let $e_i=a+(1-a)^{1\over 2} f_i (1-a)^{1\over 2}$.
Then $(e_i)$ is an approximate identity and $\overline p\leq e_i\leq 1$.
Now suppose $h$ is in $\tilde A_{sa}$ and $h\geq p$.
Then $e_i he_i\geq e_i pe_i=p$.
Therefore $K^2 e_i he_i\geq\overline p$.
Taking $\sigma$-strong limits in $A^{**}$, we see that $K^2 h\geq\overline p$.
Then (b) implies $p$ is $K$-quasi-regular.

(a$'$) and (b$'$) follow by letting $K\to 1^+$.
\enddemo

\proclaim{Theorem 4.5}If $\alpha(p)<\infty$, then $p$ is cone-regular if and only if $pbp\geq p\Rightarrow \overline p b\overline p\geq\overline p$, $\forall b\in A_{sa}$.
\endproclaim

\demo{Proof}Assume $p$ is cone-regular and $pbp\geq p$.
Then
$$
\gathered
p(1-b) p\leq 0,\\
p(e_i-b) p\leq 0,\ \forall i,\text{ where }(e_i)\text{is an approximate identity}\\
\overline p(e_i-b)\overline p\leq 0,\ \forall i,\text{ by }4.1(d),\\
\overline p(1-b)\overline p\leq 0,\text{ by taking the }\sigma-\text{strong limit and,}\\
\overline p b\overline p\geq \overline p.
\endgathered
$$

Next assume $pap\geq p$ (possible since $\alpha(p) <\infty$),\newline
$pbp\geq p\Rightarrow\overline p b\overline p\geq\overline p$, $\forall b\in A_{sa}$, and $pxp\leq 0$, $x\in A_{sa}$.
Then
$$
\gathered
p(a-tx)p\geq p,\ \forall t>0,\\
\overline p(a-tx)\overline p\geq\overline p,\ \forall t>0,\\
\overline p x\overline p\leq t^{-1} \overline p(a-1)\overline p,\ \forall t>0,\text{ and hence}\\
\overline p x\overline p\leq 0,\text{ by taking the limit as }t\to\infty.
\endgathered
$$
By 4.1(d), the above shows $p$ is cone-regular.
\enddemo

\proclaim{Corollary 4.6}If $p$ is cone-regular, then $\alpha(p)=\alpha(\overline p)$.
\endproclaim

\proclaim{Theorem 4.7}Assume $p$ is cone-regular.\newline
(a)\ If $\alpha(p)=1$, then $p$ is regular.\newline
(b)\ If $\alpha(p)=s<\infty$, then $p$ satisfies $(R_8(s))$.
A fortiori $p$ satisfies $(R_9(s))$ and $p$ is $s^{1\over 2}$-quasi-regular.
\endproclaim

\demo{Proof}(b)\ If $\varphi\in F(\overline p)$, there is a net $(\varphi_i)$ in $C(p)$ such that $\varphi_i\overset{w^*}\to\rightarrow \varphi$.
Let $\psi_i={\varphi_i\over \|\varphi_i\|}$, an element of $F(p)\cap S(A)$.
If $\|\varphi_i\| \rightarrow \infty $, then $\psi_i \overset{w^*}\to\rightarrow 0$, in contradiction to $\alpha(p)<\infty$ and 2.10.
Thus, passing to a subnet if necessary, we may assume $\|\varphi_i\|\to t<\infty$.
Then $t^{-1}\varphi\in\overline S(p)$.
By 2.9.
$s\geq t$.
By definition, this shows $(R_8(s))$.

(a)\ follows from (b) if we let $s\to 1^+$.

If $A$ is non-unital, we can identify $\tilde A^{**}$ with $A^{**}\oplus\bC$.
Then any projection in $A^{**}$ can also be regarded as an element of $\tilde A^{**}$.
Of course, the next result is also true, trivially, if $A$ is unital.
\enddemo

\proclaim{Theorem 4.8}If $p$ is a projection in $A^{**}$, then $p$ is regular in $\tilde A^{**}$ if and only if $p$ is regular in $A^{**}$ and $\alpha(p)$, computed relative to $A$, is 1 or $\infty$.
\endproclaim

\demo{Proof}We first reduce to the case where $p$ is closed, and hence regular, in $A^{**}$.
Let $p_1$ and $p_2$ be the closures of $p$ in $A^{**}$ and $\tilde A^{**}$ respectively.
It is easy to see (cf.~[7, 3.54]) that $p_1\leq p_2$.
If $p$ is regular in $A^{**}$, then by the weak$^*$ lower semicontinuity of norm, $[F(p)\cap S(A)]^- \supset F(p_1)\cap S(A)$.
Since the weak$^*$ topologies of $A^*$ and $\tilde A^*$ agree on $S(A)$ this shows that $F(p)$ and $F(p_1)$ have the same weak$^*$ closures in $\tilde A^*$.
Thus $p$ is regular in $\tilde A^{**}$ if and only if $p_1$ is.
Also $\alpha(p)=\alpha(p_1)$ by 4.6.
Now assume $p$ is regular in $\tilde A^{**}$.
Then the $\tilde A^*$-closure of $F(p)\cap S(A)$ includes $F(p_2)\cap S(\tilde A)$, and hence it includes $F(p_1)\cap S(A)$.
Again since the two weak$^*$ topologies agree on $S(A)$, the $A^*$-closure of $F(p)\cap S(A)$, includes $F(p_1)\cap S(A)$; and hence $p$ is regular in $A^{**}$.

From now on we assume $p$ closed in $A^{**}$ and let $\overline p$ denote $p_2$.
If $\alpha(p)=1$, $p$ is compact in $A^{**}$.
This implies by results of Akemann [4] (cf.~also [7, 2.47]) that $p$ is closed in $\tilde A^{**}$, and hence regular in $\tilde A^{**}$.
If $\alpha(p) > 1$, then $p$ is not compact and, by the results cited above, $p$ is not closed in $\tilde A^{**}$.
Thus $\overline p=p\oplus 1$ in $A^{**}\oplus\bC$.
If $p$ is regular in $\tilde A^{**}$, then $0\oplus 1$ is in the weak$^*$ closure of $F(p)\cap S(A)$, where $\tilde A^*$ is identified with $A^*\oplus\bC$.
This implies $0\in\overline S(p)$, and by 2.10, $\alpha(p)=\infty$.
If $\alpha(p)=\infty$, then by 2.10, 0 is in $\overline S(p)$.
Since $\overline S(p)$ is convex, this shows $\overline S(p)=F(p)$.
Now $S(\tilde A)$, with its weak$^*$ topology, can be identified with $Q(A)$, with its weak$^*$ topology.
The map is \newline
$\varphi\leftrightarrow \varphi\oplus (1-\|\varphi\|)$.
Thus $\overline S(p)=F(p)$ implies that the $\tilde A^*$-closure of $F(p)\cap S(A)$ is $F(\overline p)\cap S(\tilde A)$.
This shows that $p$ is regular in $\tilde A^{**}$.
\enddemo

\proclaim{Theorem 4.9}If $p$ is $K$-quasi-regular for $K<\sqrt{2}$, then $p$ satisfies $(R_9({K^2\over 2-K^2}))$.
\endproclaim

\demo{Proof}Use 4.1(b) to interpret $(R_9(\cdot))$.
Assume $a\in A_{sa}$, $\|pap\|=1$, and $\|\overline p a\overline p\|=s$.
By 3.3 or 3.4 of [7], there is $b$ in $A_{sa}$ such that $\|b\|=s$ and $\overline p b\overline p=\overline p a\overline p$.
Therefore also $pbp=pap$.
Assume that the top point in $\sigma(\overline p a\overline p)$ is $s$ (otherwise replace $a$ by $-a$ and $b$ by $-b$).
Since $b+s\in\overline{A_+^m}, \|\overline p(b+s)\overline p\|\leq K^2 \|p(b+s)p\|$, by 4.2(b).
Thus $2s\leq K^2(1+s)$, and hence $s\leq {K^2\over 2-K^2}$.
\enddemo

\medskip\noindent
{\bf Examples 4.10}.
(a)\ First of all, we promised at the beginning of Section 3 to say what effect (ordinary) regularity has on $(\alpha(p),\alpha(\overline p))$.
By 4.6, if $p$ is regular, $\alpha(p)=\alpha(\overline p)$.
Thus we consider regularity only in connection with 3.3, 3.4, and 3.5.
By a result of Tomita [21], every central projection is regular.
Thus our ``examples'' in 3.3 and 3.4 are regular.
In 3.5 we gave two examples with $\alpha(p)=\alpha(\overline p)=s$, $1< s<\infty$.
One example was closed, and hence regular.
The other example was open but not (cf.~4.10(b), below) regular.
It is easy to modify this example and obtain a regular and even $k$-regular, $\forall k$, open projection with $\alpha(p)=\alpha(\overline p)=s$.
Let $A$ and $v_n$ be as in 3.5, and define $p$ by $p_\infty=0$, $p_{2n}=v_n\times v_n$, $p_{2n-1}=e_1\times e_1$.
As before, $\overline p$ differs from $p$ only in that $(\overline p)_\infty=e_1\times e_1$.
It is easy to check that $p$ has the properties claimed.
For regularity, we need that $\|a_\infty(\overline p)_\infty\|\leq \sup_n \|a_n p_n\|$, $\forall a\in A$.
This follows from \newline
$\|a_\infty (e_1\times e_1)\|\leq\sup_n \|a_{2n-1}(e_1\times e_1)\|$, which is true because $a_{2n-1}\to a_\infty$ in norm.
$k$-regularity is proved similarly.

(b)\ In this example $p$ is cone-regular but not regular, and $\alpha(p)=s$, $1< s < \infty$.
Also this example shows that the estimates in 4.7(b) for the constants in $(R_i(\cdot))$ are sharp.
This example is exactly the open example given in 3.5.
If $pap\leq 0$, $a$ in $A_{sa}$, then $(a_n v_n,v_n)\leq 0$.
Since $a_n\to a_\infty$ in norm, $a_\infty\in\cK$, and $v_n\overset w\to\rightarrow s^{-{1\over 2}}e_1$ we can take a limit and obtain $s^{-1}(a_\infty e_1,e_1)\leq 0$.
Therefore $\overline p a\overline p\leq 0$.
By 4.1(d) $p$ is cone-regular.
If we define $a$ in $A$ by $a_n=a_\infty=e_1\times e_1$, then $\|a\overline p\|=1$.
$\|ap\|=\sup|(e_1,v_n)|=s^{-{1\over 2}}$.
This shows that $p$ is at best $s^{1\over 2}$-quasi-regular, precisely in accordance with 4.7(b).
It follows \underbar{ a fortiori} that the $(R_8(\cdot))$ and $(R_9(\cdot))$ estimates given in 4.7(b) are also sharp.

(c)\ In this example $p$ is cone-regular and open but not regular and not even quasi-regular.
By 4.7(b), $\alpha(p)=\infty$.
Let $A=c\otimes\cK$.
Let $\cV_m=\{u\in H\colon \|u\|=m^{-1}$ and $(u,e_k)=0$, $\forall k> m\}$, $m=1,2,\ldots$
Let $v_n$ be a sequence of unit vectors in $H$ such that:

(i)\ $\forall n$, $\exists m$ such that $v_n=u_n+(1-m^{-2})^{1\over 2} e_k$, where $u_n\in\cV_m$ and $k>\max (m,n)$, and

(ii)\ $\{u_n\}$ contains a dense subset of each $\cV_m$.

\noindent
Define an open projection $p$ in $A^{**}$ by $p_\infty=0$ and $p_n=v_n\times v_n$.
For each $u$ in $\cV_m$, $u\times u$ is a weak cluster point of $(p_n)$.
Therefore $(\overline p)_\infty=1$ (and, as always, $(\overline p)_n=p_n$ for $n<\infty$).
If $a\in A_{sa}$ and $pap\leq 0$, then $(a_n v_n,v_n)\leq 0$, $\forall n$.
As in (b), it follows that $(a_\infty u,u)\leq 0$, $\forall u\in\cV_m$.
Therefore $a_\infty\leq 0$ and $\overline pa\overline p\leq 0$.
Therefore $p$ is cone-regular.
Now let $x=e_m\times e_m$.
Then $xu=0$ if $u\in\cV_{m'}$, $m'< m$, $\|xe_k\|\to 0$ as $k\to\infty$, and $\|x\|=1$.
Thus $\limsup\|xv_n\|\leq m^{-1}$.
Now define $a_N$ in $A$ by $(a_N)_n=\cases 0,&n\leq N\\ x,&N+1\leq n\leq\infty\endcases$.
Then $\|a_N\overline p\|=1$ and $\limsup_{N\to\infty} \|a_N p\|\leq m^{-1}$.
Since $m$ is arbitrary, $p$ is not quasi-regular.

(d)\ In this example $A$ is unital and $p$ is $\sqrt{2}$-quasi-regular but not 0-regular.
This shows that the constant $\sqrt{2}$ in 4.9 is sharp.
Let $A=c\otimes M_2$.
Define an open projection $p$ in $A^{**}$ by $p_\infty=0$, $p_{2k-1}=e_1\times e_1$, and $p_{2k}=e_2\times e_2$.
Then $(\overline p)_\infty=1$.
For $a$ in $A$, $\|ap\|=\max (\sup_k \|a_{2k-1} e_1\|$, $\sup_k \|a_{2k}e_2\|)$.
Since $a_n\to a_\infty$ in norm,\newline
$\|ap\|\geq\max(\|a_\infty e_1\|$, $\|a_\infty e_2\|)$.
We conclude easily, for example by looking at the Hilbert-Schmidt norm, that $\|a_\infty\|\leq\sqrt{2}\|ap\|$.
Therefore $\|a\overline p\|\leq\sqrt{2}\|ap\|$.
If we define $b$ in $A_{sa}$ by $b_n=b_\infty=\pmatrix 0&1\\ 1&0\endpmatrix$, then $pbp=0$ and $\overline p b\overline p\neq 0$.
Therefore $p$ is not 0-regular.

\proclaim{Theorem 4.11}If $\alpha(p)=s$ and $p$ is $K$-quasi-regular with $K^2< s/(s-1)$, then $\alpha(\overline p)\leq s/[s-K^2(s-1)]<\infty$.
In particular, if $s=1$, then $\alpha(\overline p)=1$.
\endproclaim

\demo{Proof}Choose $s_1>s$ such that $K^2< s_1/s_1-1$, and choose $a$ in $A_+$ such that $pap\geq p$ and $\|a\|\leq s_1$.
Then $p(1-a) p\leq 0$ and hence $p(s_1-a) p\leq (s_1-1)p$.\newline
Since $s_1-a\in\overline{A_+^m}$, 4.2(b) implies $\overline p(s_1-a)\overline p\leq K^2 (s_1-1)\overline p$.
Therefore \newline
$\overline p a\overline p\geq [s_1-K^2 (s_1-1)]\overline p$.
Since $s_1 -K^2 (s_1-1) > 0$, this shows \newline
$\alpha(\overline p)\leq s_1/s_1 -K^2 (s_1-1)$.
Now let $s_1\to s^+$.
\enddemo

\proclaim{Lemma 4.12}Let $H$ be a Hilbert space, $e_1$ a unit vector in $H$, and $Q$ the projection with range $\{e_1\}^\perp$.
Let $0< t<1$, $W_1=\{u\in H\colon \|u\|=1$ and $\|Qu\|\leq t\}$, and $W=\{u\in H\colon \|u\|\leq 1$ and $\|Qu\|\leq t\}$.
Then $W$ is a balanced convex set and is the closed convex hull of $W_1$.
\endproclaim

\demo{Proof (sketch)}$W$ is weakly compact and $W_1$ is the set of extreme points.
\enddemo

\noindent
{\bf Examples 4.13}.\ (a)\ In this example $\alpha(p)=s$, $1<s<\infty$, $p$ is open and $K$-quasi-regular with $K^2={s\over s-1}$, and $\alpha(\overline p)=\infty$.
This shows that the estimate on $K$ in 4.11 is sharp.
This also shows that $K$-quasi-regularity does not imply cone-regularity for $ K>1$.

Let $A=c\otimes\cK$, $t=(1-s^{-1})^{1\over 2}=K^{-1}$, and let $(v_n)$ be a dense sequence in $W_1$.
Define an open projection $p$ in $A^{**}$ by $p_\infty=0$ and $p_n=v_n\times v_n$.
As in 3.6, we see that $(\overline p)_\infty=1$, so that $\alpha(\overline p)=\infty$.
If $a$ in $A$ is defined by $a_n=a_\infty=e_1\times e_1$, then $pap\geq s^{-1} p$.
Therefore $\alpha(p)\leq s$.
If $b\in A$ and $\|bp\|\leq 1$, then $\|b_n v_n\|\leq 1$, $\forall n$.
Therefore $\|b_\infty u\|\leq 1$, $\forall u\in W_1$.
By 4.12, $\|b_\infty u\|\leq 1$, $\forall u\in W$, and hence $\|b_\infty\| < t^{-1}=K$.
We have shown $\|b\overline p\|\leq K$, and thus $p$ is $K$-quasi-regular.
4.11 now shows that $\alpha(p)=s$.

(b)\ If we want a unital example where $K$-quasi-regularity does not imply (cone-) regularity, or better, if we want $p$ to be $K$-quasi-regular and not $K'$-quasi-regular for any $K'< K$, we can use the same construction as in (a) for $A=c\otimes M_2$.
Thus now the $H$ of 4.12 is two dimensional and $\alpha(p)=\alpha(\overline p)=1$.
Since $W$ contains the ball of radius $t$ but no larger balls, the separation theorem shows that for any $t' > t$ we can find a linear functional 
$h$ on $H$ such that $|h(u)|\leq 1$, $\forall u\in W$, and $|h(u_0)| > 1$ for some $u_0$ with $\|u_0\|=t'$.
Define $a$ in $A$ by $a_n=a_\infty=e_1\times y$, where $y$ in $H$ is such that $h(\cdot)=(\cdot,y)$.
Then $\|ap\|\leq 1$ and $\|a\overline p\|\geq (t')^{-1}$.

(c)\ This example shows that the estimate on $\alpha(\overline p)$ in 4.11 is sharp when $K^2< {s\over s-1}$.
The construction is similar to 3.9, but unfortunately it must be a bit more complicated if we want $p$ to be open.
Thus we assume given $s$ and $t$ such that $1< s<t<\infty$ and let $s'$ be as in 3.9.
Let $K^2={s'\over s'-1}$.
The reader can compute that $t={s\over s-K^2(s-1)}$.
(If only $s$ is given, $t$ can be chosen so that $K$ takes any arbitrary value in $(1,({s\over s-1})^{1\over 2}$).

Thus we perform the construction of (a) with $s'$ instead of $s$.
Let $A_0$ and $p_0$ be the algebra and projection produced by this, and let $A_1,e',A,e$, and $q$ have the same meaning as in 3.9 (and 3.8).
(We now have $p_0 q p_0\geq (s')^{-1} p_0$ instead of equality.)

For each $n$ choose a unit vector $z_n=x_n\oplus y_n$ in $H\oplus H$ such that:
\item{(i)}$(x_n,e_1)=(y_n,e_1)=0$,
\item{(ii)}$(x_n,v_n)=0$,
\item{(iii)}$e' z_n=z_n$,
\item{(iv)}$z_n\overset w\to\rightarrow 0$ as $n\to\infty$.

\noindent
Let $p$ be the open projection in $A_1^{**}$ defined by $p_\infty=0$ and\newline
$p_n=(v_n\oplus 0)\times (v_n\oplus 0)+z_n\times z_n$.
As before, $p$ is also regarded as an open projection in $A^{**}$, the closure, $\overline p$, of $p$ in $A_1^{**}$ is the same as $p$ except that $(\overline p)_\infty=\pmatrix 1&0\\ 0&0\endpmatrix$, and the closure of $p$ in $A^{**}$ is $\overline p\oplus 1$.

Let $a$ have the same meaning as in 3.9 ($a$ is a specific element of norm 1 in $A_+$).
Then using (i), (ii), (iii), we see, similarly to 3.9, that $pap\geq s^{-1} p$ and\newline
$(\overline p\oplus 1) e(\overline p\oplus 1)\geq t^{-1} (\overline p\oplus 1)$.
Also, the proofs in 3.8 and 3.9 that $\alpha(p)\geq s$ and $\alpha(\overline p\oplus 1)\geq t$ still apply, since $p\geq \pmatrix p_0&0\\ 0&0\endpmatrix\oplus 0$.

It remains to show that $p$ is $K$-quasi-regular.
If $b\in A$, since the $A_1^{**}$-component of $b$ is in $M(A_1)$, the calculations of (a) (cf.~4.3(b)) show that\newline
$\|b\overline p\|\leq ({s'\over s'-1})^{1\over 2} \|bp\|$.
(Note that to show this, we need only estimate $\|b_\infty(\overline p)_\infty\|$, and we do not need to consider the $z_n$'s).
If $b=\lambda e+x$, $x\in A_1$, then by (iii) and (iv) $\|bp\|\geq |\lambda|$.
This concludes the proof.

If we drop the openness requirement, we can get an easier example:\newline
Let $p=\pmatrix p_0&0\\ 0&0\endpmatrix\oplus 1$.
In this case $p$ is abelian, as in many of our earlier examples.

The gist of what we have done so far is that in general cone-regularity and quasi-regularity are independent of one another, and that both have significant relations with near relative compactness.
A special case of our results is that if $p$ is either cone-regular or quasi-regular, then condition (2) of Section 1 implies $p$ is relatively compact (briefly, $\alpha(p)=1\Rightarrow \alpha(\overline p)=1$).

In Section 6 we will consider situations in which $\alpha(p_1\vee p_2)$ can be bounded in terms of $\alpha(p_1)$ and $\alpha(p_2)$.
For the question of regularity of $p_1\vee p_2$, we will consider only the special case where $\overline p_1 \overline p_2=0$, so that $p_1\vee p_2=p_1+p_2$ and $(p_1\vee p_2)^- =\overline p_1+\overline p_2$.
Also we consider only the hypothesis that $p_1$ and $p_2$ are regular in the ordinary sense, except when generalizations are easy.
The conclusions available even from these seemingly strong hypotheses are not strong.
Of course, we do not expect to be able to prove $p_1+p_2$ is regular---otherwise there would be no purpose for the concept of $k$-regularity.

\proclaim{Proposition 4.14}(a)\ If $p_i$ is $K_i$-quasi-regular for $i=1,2$, and if $\overline p_1\overline p_2=0$, then $p_1+p_2$ is $(K_1^2+K_2^2)^{1\over 2}$-quasi-regular.

(b)\ If $p$ is 0-regular, then Diag$(p,\ldots,p)$ is 0-regular in $(A\otimes M_k)^{**}$.
\endproclaim

\demo{Proof}(a)\ By [1], $\overline p_1+\overline p_2$ is closed, so that $(p_1+p_2)^-=\overline p_1+\overline p_2$.
If $\|a(p_1+p_2)\|\leq 1$ for $a$ in $A$, then $\|a\overline p_i\|\leq K_i$.
Thus $a(\overline p_1+\overline p_2)a^*\leq a\overline p_1 a^* + a\overline p_2 a^*\leq K_1^2+K_2^2$,\newline
and $\|a(\overline p_a+\overline p_2)\|\leq (K_1^2+K_2^2)^{1\over 2}$.

(b)\ Since Diag$(p,\ldots,p)^- =$ Diag$(\overline p,\ldots,\overline p)$, this follows immediately from 4.1(c).
\enddemo

\noindent
{\bf Examples 4.15}.\ (a)\ Let $\pi\colon M_2\to B(H)$ be a unital $*$-representation, so that $\pi$ induces a faithful homomorphism from $M_2$ to the Calkin algebra, $B(H)/\cK$.
Let $A$ be the extension of $\cK\oplus\cK$ by $M_2$ induced by $\pi\oplus\pi\colon M_2\to B(H)\oplus B(H)$.
(cf.~[10].)
Then $A^{**}$ can be identified with $B(H)\oplus B(H)\oplus M_2$.
Let $\{e_{ij}\colon i,j=1,2\}$ be a system of matrix units for $M_2$, and define $p_1=\pi(e_{11})\oplus 0\oplus 0$, $p_2=0\oplus\pi(e_{22})\oplus 0$.
It is easy to check that $\overline p_1=\pi(e_{11})\oplus 0\oplus e_{11},\overline p_2=0\oplus\pi(e_{22})\oplus e_{22}$, and $p_i$ is $k$-regular, $\forall k$, and open.
If $a_0$ is $\pmatrix 0&1\\ 1&0\endpmatrix$ in $M_2$ and $a=\pi(a_0)\oplus \pi(a_0)\oplus a_0$, then $a\in A$, $(p_1+p_2)a(p_1+p_2)=0$, and $(\overline p_1+\overline p_2) a(\overline p_1+\overline p_2)\neq 0$.
Thus $p_1+p_2$ is not even 0-regular.
Thus $p_1+p_2$ cannot be better than $\sqrt{2}$-quasi-regular by 4.9, so that 4.14(a) is sharp, at least in the special case $K_1=K_2=1$.
Note that $A$ is unital.

(b)\ The fact that for every $k>1$ there is a projection which is $(k-1)$-regular but not $k$-regular is surely due to the inventor of $k$-regularity.
For completeness, we write down a natural example, but the proof that it is correct is left to the reader.
Let $A=c\otimes M_k$, a unital algebra, and let $(q_n)$ be a sequence dense in the set of rank $k-1$ projections in $M_k$.
Define a $(k-1)$-regular open projection $p$ in $A^{**}$ by $p_\infty=0$ and $p_n=q_n$.
Then $(\overline p)_\infty=1$ and $p$ is not $k$-regular.

\proclaim{Theorem 4.16}If $\overline p_1\overline p_2=0$, then $p_1+p_2$ is regular if and only if:
\item{(i)}$p_1$ and $p_2$ are regular, and
\item{(ii)}$\{f\in A^*\colon \|f\|\leq 1$ and $f(a)=f(p_1 ap_2),\forall a\}$ is weak$^*$ dense in $\{f\in A^*\colon \|f\|\leq 1$ and $f(a)=f(\overline p_1 a\overline p_2),\ \forall a\}$.
\endproclaim

\demo{Proof}We first assume (i) and (ii) and prove $(R_2), F(p_1+p_2)^- =F(\overline p_1+\overline p_2)$.
Since $F(p_1+p_2)^-$ is convex, it is enough to show it contains each pure state $\varphi$ in $F(\overline p_1+\overline p_2)$.
Let $\pi=\pi_\varphi$.
Then $\varphi=(\pi(\cdot)v,v)$ for some unit vector $v$ in $H_\pi$, $v=v_1+v_2$, where $v_i=\pi^{**} (\overline p_i)v$.
If $v_1$ or $v_2$ is 0, then $\varphi$ is in $F(\overline p_2)$ or $F(\overline p_1)$, and by (i) $\varphi$ is in $F(p_1+p_2)^-$.
Thus assume $v_1,v_2\neq 0$, and let $u_i={v_i\over \|v_i\|}$.
If we define $f$ in $A^*$ by $f(\cdot)=(\pi(\cdot) u_2,u_1)$, then $\|f\|=1$ since $\pi$ is irreducible (use [15]), and clearly $f(\cdot)=f(\overline p_1\cdot\overline p_2)$.

Consider $A^{**}$ as a subset of $B(H)$, via the universal representation of $A$.
(This is just a matter of convenience in ``bookkeeping''.)
By (ii), we can find a net $(f_i)$ such that $\|f_i\|=1$, $f_i(\cdot)=f_i(p_1\cdot p_2)$, and $f_i\overset{w^*}\to\longrightarrow f$.
Write $f_i(\cdot)=(\cdot w_i^2,w_i^1)$, where $w_i^j$ is a unit vector in $p_j H$, $j=1,2$.
We obtain this representation of $f_i$ from the polar decomposition ([20]).
Note that $|f_i|=(\cdot w_i^2,w_i^2),|f_i^*|=(\cdot w_i^1,w_i^1)$, $|f|=(\pi(\cdot) u_2,u_2)$, and $|f^*|=(\pi(\cdot) u_1,u_1)$.
Since $f_i\overset{w^*}\to\longrightarrow f$ and $\|f_i\|\to \|f\|$, it follows from [13] that $|f_i|\overset{w^*}\to\longrightarrow |f|$ and $|f_i^*|\overset{w^*}\to\longrightarrow |f^*|$.
Now let $w_i=\|v_1\|w_i^1+\|v_2\|w_i^2$ and $\varphi_i=(\cdot w_i,w_i)$.
Then $\varphi_i\in F(p_1+p_2)$ and $\varphi_i\overset{w^*}\to\longrightarrow \varphi$.

Next assume $p_1+p_2$ is regular.
Every element $f$ of $V(\overline p_1+\overline p_2)$ can be written uniquely as $\sum\limits^2_{j,k=1} f^{jk}$, where $f^{jk}(\cdot)=f^{jk}(\overline p_j\cdot\overline p_k)$.
Also a bounded net $(f_i)$ in $V(\overline p_1+\overline p_2)$ converges weak$^*$ to $f$ if and only if $f_i^{jk}\overset{w^*}\to\longrightarrow f^{jk}$, $\forall j,k$.
Of course, $F(\overline p_1+\overline p_2)\subset V(\overline p_1+\overline p_2)$, so that the above applies.
If $\varphi\in F(\overline p_1)\cap S(A)$, then by $(R_2)$, there is a net $(\varphi_i)$ in $F(p_1+p_2)$ such that $\varphi_i\overset{w^*}\to\longrightarrow \varphi$.
Then also $\varphi_i^{11}\overset{w^*}\to\longrightarrow\varphi$.
Since $\varphi_i^{11} (\cdot)=\varphi_i (\overline p_1\cdot\overline p_1)$ and $\varphi_i\in F(p_1+p_2)$, $\varphi_i^{11}\in F(p_1)$.
Therefore $p_1$, and similarly $p_2$, are regular.

Now assume $\|f\|=1$ and $f(\cdot)=f(\overline p_1\cdot\overline p_2)$.
Write $f(\cdot)=(\cdot u_2,u_1)$, where $u_j$ is a unit vector in $\overline p_j H$.
Let $v=2^{-{1\over 2}} u_1+2^{-{1\over 2}}u_2$, and $\varphi=(\cdot v,v)$.
Then $\varphi\in F(\overline p_1+\overline p_2)\cap S(A)$.
By $(R_2)$ there is a net $(\varphi_i)$ in $F(p_1+p_2)$ such that $\varphi_i\overset{w^*}\to\longrightarrow\varphi$.
Therefore $\varphi_i^{jk}\overset{w^*}\to\longrightarrow\varphi^{jk}$.\newline
Then $2^{-1}=\|\varphi^{jj}\|\leq\liminf\|\varphi_i^{jj}\|$.
Since $\|\varphi_i^{11}\|+\|\varphi_i^{22}\|=\|\varphi_i\|\leq 1$, it follows that \newline
$\|\varphi_i^{jj}\|\to 2^{-1}$.
Now $\|\varphi_i^{12}\|\leq \|\varphi_i^{11}\|^{1\over 2} \|\varphi_i^{22}\|^{1\over 2}$, and\newline
$\varphi_i^{12}\overset{w^*}\to\longrightarrow \varphi^{12}=2^{-1} f$.
Therefore $\|\varphi_i^{12}\|\to 2^{-1}$, and ${\varphi_i^{12}\over \|\varphi_i^{12}\|}\to f$.
Finally, since $\varphi_i\in F(p_1+p_2)$ and $\varphi_i^{12}=\varphi_i(\overline p_1\cdot\overline p_2)$, $\varphi_i^{12}=\varphi_i^{12} (p_1\cdot p_2)$.
Thus (ii) is proved
\enddemo

\proclaim{Corollary 4.17}2-regular $\Leftrightarrow (R_5)$.
Explicitly, $p$ is 2-regular if and only if\newline
$\{f\in A^*\colon \|f\|\leq 1$ and $f(a)=f(pap),\ \forall a\}$ is weak$^*$ dense in\newline
$\{f\in A^*\colon \|f\|\leq 1$ and $f(a)=f(\overline pa\overline p),\ \forall a\}$. 
\endproclaim

\demo{Proof}By 4.3(f), $(R_5)\Rightarrow$ regular.
\enddemo

\remark{Remark}The boundedness hypothesis on the net $(f_i)$ in the proof of 4.16 is not really needed (to deduce $f_i^{jk}\to f^{jk}$), at least when $A$ is $\sigma$-unital.
A proof can be based on the Urysohn lemma ([7.\ 3.31]).
\endremark

\subheading{\S 5.\ Special Results for Open Projections}

\proclaim{Theorem 5.1}Suppose $p$ and $q$ are projections in $A^{**}$ such that $p$ is closed, $q$ is open, $\alpha(p)<\infty$, and $\|p-q\| < 1$.
Then $q$ is compact.
Thus $q$ is actually in $A$.
\endproclaim

\demo{Proof}Since $\|p-q\|<1$, there is $\var>0$ such that $pqp\geq\var p$.
Let $B$ be the hereditary $C^*$-subalgebra of $A$ supported by $q$ (notation:\ $B=\her(q)$), and let $(e_i)$ be an approximate identity of $B$.
Then $pe_i p\nearrow pqp$, with convergence in the $\sigma$-strong topology.
Since $p$ is closed, $\overline S(p)$ is a weak$^*$ compact subset of $F(p)$, and $\|\varphi\|\geq\alpha (p)^{-1}$ for $\varphi$ in $\overline S(p)$.
Therefore, $\forall\varphi\in\overline S(p)$, $\lim\varphi(e_i)=\lim\varphi (pe_i p)\geq\var \alpha (p)^{-1}$.
By Dini's theorem, for $i$ sufficiently large $\varphi (e_i)\geq 2^{-1}\var\alpha (p)^{-1}$, $\forall\varphi\in \overline S(p)$.
Since $\overline S(p)\supset F(p)\cap S(A)$, we have shown $pe_i p\geq\delta p$ for $i\geq i_0$, where $\delta=2^{-1}\var\alpha (p)^{-1} > 0$.
\enddemo

Now let $p_i$ be the range of projection of $e_i^{1\over 2} p$ for $i\geq i_0$.
By the proof of Theorem 2.6(b), $p_i$ is compact and $\|p_i-p\| < 1$.
Also, clearly, $p_i\leq q$.
Since also $\|q-p\| < 1$, it follows that $p_i=q$ (cf.~Lemma 7 of [9] and the remark following).

Since $q$ is both closed and open $q\in M(A)$.
It is easy to see that a projection in $M(A)$ is compact, as an element of $A^{**}$, if and only if it is in $A$.

\medskip\noindent
{\bf Example 5.2}.
There are a closed but not open projection $p$ and an open but not closed projection $q$ such that $\alpha(q) <\infty$ and $\|p-q\|<1$.
Let $A=c\otimes\cK$.
Let $v_n=2^{-{1\over 2}} e_2+2^{-{1\over 2}} e_{n+2}$ and $w_n=2^{-{1\over 2}} e_1+2^{-{1\over 2}} e_{n+3}$.
Define $p$ and $q$ by $p_\infty=q_\infty=e_1\times e_1$,\newline
$p_n=w_n\times w_n+e_{n+2}\times e_{n+2}$, and $q_n=e_1\times e_1+v_n\times v_n$.
Then $\alpha(q)=2$ and $\|p-q\|=2^{-{1\over 2}}$.  $(d_a (q,RC)={\pi\over 4}$, $d_a(p,q)={\pi\over 4}$.)

We have already mentioned that if $A$ is $\sigma$-unital and $q$ is an open projection in $A^{**}$, then $\alpha(q)<\infty$ if and only if every closed subprojection of $q$ is compact.
This results from the combination of Theorems 2 and 4 of [8], and its proof there is indirect:\ It is shown that both conditions are equivalent to $M(A,B)=B$, where $B=\her(q)$ and $M(A,B)=M(A)\cap (qA^{**}q)$.
The proof of one half of this result in [8] is actually quite easy, granted the Urysohn lemma, [7,3.31].
Nevertheless, we give a more direct proof of this half.

\proclaim{Lemma 5.3}If $A$ is $\sigma$-unital, $p$ and $q$ are projections in $A^{**}$, $p$ is closed, $q$ is open, $p\leq q$, and $\alpha(q)<\infty$, then $p$ is compact.
\endproclaim

\demo{Proof}By [7.~3.31] there is $h$ in $M(A)_{sa}$ such that $p\leq h\leq q$.
Choose $a$ in $A_{sa}$ such that $qaq\geq q$.
Then $h^{1\over 2} ah^{1\over 2}=h^{1\over 2} qaqh^{1\over 2}\geq h^{1\over 2} qh^{1\over 2}=h\geq p$.
Since $h^{1\over 2}ah^{1\over 2}\in A$, $p$ satisfies (2) of Section 1, and by Akemann [4] $p$ is compact.
\enddemo

\proclaim{Theorem 5.4}If $A$ is any $C^*$-algebra, $p$ and $q$ are projections in $A^{**}$, $p$ is closed, $q$ is open, $p\leq q$, and $\alpha(q)<\infty$, then $p$ is compact.
\endproclaim

\demo{Proof}The proof is by reduction to the separable case.
Choose $a$ in $A_{sa}$ such that $qaq\geq q$.
Let $q'=1-p$, $B=\her(q)$, and $B'=\her(q')$.
Choose $\sigma$-unital hereditary $C^*$-subalgebras $B_0$ of $B$ and $B'_0$ of $B'$ such that $a\in\her (B_0\cup B'_0)$.
(For a subset $S$ of $A$, $\her(S)$ denotes the smallest hereditary $C^*$-subalgebra including $S$.)
Let $e_0,e'_0$ be strictly positive elements of $B_0,B'_0$, and let $A_0=C^*(e_0,e'_0,a)$, a separable $C^*$-subalgebra of $A$.
If $p$ is not compact, there is a net $(\varphi_i)_{i\in D}$ in $F(p)\cap S(A)$ such that $\varphi_i\overset{w^*}\to\longrightarrow \varphi$ and $\|\varphi\|<1$.
The construction of the desired separable $C^*$-subalgebra of $A$ proceeds from here by recursion.
\enddemo

Step 1.
Since $A_0$ is separable, we can choose $i_1,i_2,\ldots$ such that\newline
$\varphi_{i_n}|_{A_0}\overset{w^*}\to\longrightarrow \varphi|_{A_0}$.
Since each $\|\varphi_i|B\|=1$ (because $\varphi_i\in F(p)\cap S(A)$ and $p\leq q$), we can find a countable subset $E_1$ of $\{b\in B\colon 0\leq b\leq 1\}$ such that $\sup\varphi_{i_n}|_{E_1}=1$ for each $n$.
Then, by [7, 3.30], we can find open projections $q_1,q'_1$ such that $q_1\leq q$, $q'_1\leq q', q_1 q'_1=q'_1 q_1$, $\her(q_1)$ and $\her(q'_1)$ are $\sigma$-unital, and $E_1\cup \{e_0\}\subset\her(q_1)$, $e'_0\in\her(q'_1)$.
Let $e_1,e'_1$ be strictly positive elements of $\her(q_1)$, $\her(q'_1)$, and let $A_1=C^*(A_0,E_1,e_1,e'_1)$.
Note that since $A_0\subset\her(e_0,e'_0)$, $A_1\subset\her(q_1\vee q'_1)$.

Step 2 is done the same way as step 1, starting with $A_1,e_1,e'_1$ instead of $A_0,e_0,e'_0$.
(The sequence $(\varphi_{i_n})$ constructed in step 2 might be disjoint from the one in step 1.
They simply are both sequences, not \underbar{sub}sequences, constructed from the elements of the original net $(\varphi_i)$.)

The process is continued recursively, and we get increasing sequences $(q_n),(q'_n)$ of open projections and an increasing sequence $(A_n)$ of separable subalgebras.
Let $q_\infty=\lim q_n$, $q'_\infty=\lim q'_n$, and $A_\infty=C^*(\bigcup\limits_n A_n)$.
Then $q_\infty,q'_\infty$ are open, $\her(q_\infty)$, $\her(q'_\infty)$ are $\sigma$-unital, $q_\infty q'_\infty=q'_\infty q_\infty$, and $A_\infty$ is separable.
Since each $q_n$ is the range projection of $e_n$, and $e_n\in A_\infty$, $q_\infty\in A_\infty^{**}\subset A^{**}$.
Similarly, $q'_\infty\in A_\infty^{**}$.
Since $A_n\subset\her(q_n \vee q'_n)$, $q_\infty\vee q'_\infty$ is the identity of $A_\infty^{**}$.
Now $q_\infty$ and $q'_\infty$ are open as elements of $A_\infty^{**}$ ([7, 2.14(a)]), and hence if $p_\infty=(q_\infty\vee q'_\infty)-q'_\infty,p_\infty$ is a closed projection in $A_\infty^{**}$.
Since each $\varphi_i$ is supported by $p,\varphi_i|_{A_\infty}$ is supported by $p_\infty$.
Let $T=\{\varphi_i|_{A_\infty}\colon i\in D$ and $\|\varphi_i|_{A_\infty}\|=1\}$.
$T$ contains all the $\varphi_{i_n}$'s constructed in all the steps.
Thus $\varphi|_{A_\infty}$ is in the weak$^*$ closure of $T$, and $\|\varphi|_{A_\infty}\|<1$.
This shows that $p_\infty$ is not compact in $A_\infty^{**}$.
Now $a\in A_\infty$ and $qaq\geq q$ implies $q_\infty aq_\infty\geq q_\infty$.
Thus all the hypotheses of Lemma 5.3 are satisfied by $A_\infty,p_\infty,q_\infty$ and we have a contradiction.
Thus $p$ is compact after all.

\medskip\noindent
{\bf Example 5.5}.
We give a commutative counterexample to the converse of 5.4.
Of course, the example must be non-$\sigma$-unital.
Let $X$ be an ordered set with the order type of the first uncountable ordinal, endowed with the order topology.
$X$ is locally compact Hausdorff, and we let $A=C_0(X)$.
Let $U$ be the open set consisting of all isolated points of $X$ (non-limit ordinals), and let $q$ be the corresponding open projection in $A^{**}$.
Since $U$ is cofinal in $X,\alpha(q)=\infty$.
We claim that any closed subprojection, $p$, of $q$ is compact.
In fact $p$ corresponds to a closed subset $F$ of $X$ such that $F\subset U$.
$X\backslash U$ is a closed cofinal set.
Any two closed cofinal subsets of $X$ have a non-empty intersection ([16]).
Therefore $F$ is not cofinal, and $F$ and $p$ are compact.

A result of Akemann [2] states that if $a\in A$, $p$ is a closed projection in $A^{**}$, and $\|ap\|\leq 1$, then for any $\var>0$ there is an open projection $q$ such that $q\geq p$ and $\|aq\| < 1+\var$.
This is appealing from the point of view of non-commutative topology, and thus it is natural to consider similar questions.
More discussion is given after the next theorem, but we have no particular applications in mind.

\proclaim{Theorem 5.6}Assume $p$ is a closed projection in $A^{**}$, $a\in A$ and $0<\var<1$.
\item{(a)}If $a^*=a$ and $pap\leq 0$, then there is an open projection $q$ such that $q\geq p$ and $qaq\leq\var q$.
\item{(b)}If $a^*=a$ and $pap\leq p$, then there is an open projection $q$ such that $q\geq p$ and $qaq\leq (1+\var)q$.
\item{(c)}There is an open projection $q$ such that $q\geq p$ and $\|qaq\| < \|pap\|+\var$.
\item{(d)}If $a^*=a$ and $pap\geq p$, then there is an open projection $q$ such that $q\geq p$ and $qaq\geq (1-\var)q$ if and only if $p$ is compact.
\endproclaim

\demo{Proof}Let $L$ and $R$ be the closed left and right ideals of $A$ corresponding to $p$ $(L=A\her (1-p),\ R=\her(1-p)A)$.
By a result of Combes [11], $L+R$ is closed, and $L+R=\{a\in A\colon pap=0\}$.
Let $(e_i)$ be an approximate identity of $\her(1-p)$.
\item{(c)}It is known from [6] that $\|pap\|=\|a+L+R\|$ in $A/L+R$.
Thus we can find $l\in L$ and $r\in R$ such that $\|a+l+r\| < \|pap\| + {\var\over 3}$.
Then there is an $i$ such that\newline
$\|l(1-e_i)\|$, $\|(1-e_i)r\| < {\var\over 3}$.
If we let $q=E_{[0,\delta)}(e_i)$ for $\delta$ sufficiently small, then\newline
$\|lq\|$, $\|qr\| < {\var\over 3}$, and hence $\|qaq\| < \|pap\|+\var$.
\enddemo

(a) and (b) Now $a^*=a$ and we will, possibly unnecessarily, use [7, 3.4]:\ If $I$ is a closed interval that contains $\sigma(pap)\cup \{0\}$, then there is $b$ in $A_{sa}$ such that $\sigma(b)\subset I$ and $pbp=pap$.
In case (a), $I=[*,0]$ and in case (b), $I=[*,1]$.
Then, similarly to case (c), $a=b+l+r$ where $b\leq 0$ or $b\leq 1$, and we need only choose $q$ so that $\|q(l+r)q\|\leq\var$.

(d)\ If $qaq\geq (1-\var)q$, then $\alpha(q)<\infty$ and 5.4 implies $p$ is compact.
If $p$ is compact, then $p$ is closed in $\tilde A^{**}$, and we can apply (a) to $1-a$ and $\tilde A$.

The inequalities considered in (a), (b), (c), and [2] correspond to different kinds of regularity, interpreted via polars as in Section 4.
(a) relates to $(R_3)$, cone-regularity, (b) relates to $(R_2)$, regularity, (c) relates to $(R_5)$ and special cases of (c) relate to $(R_6)$ and $(R_4)$, and [2] relates to $(R_1)$, regularity.
If, for example, $p$ is a cone-regular projection and $pap\leq 0$, then $\overline p a\overline p\leq 0$ and (a) can be applied to $\overline p$.

Now 4.2 and 4.3 state that if $p$ is regular, the equality $\|xp\|=\|x\overline p\|$ is valid not only for $x$ in $A$ but also for $x$ in $QM(A)$, in particular for $x$ in $\tilde A$.
It might be hoped then that the following is true for a closed projection $p$ in $A^{**}$:

\noindent
(4)\ If $x\in\tilde A$ and $\|xp\|\leq 1$, then $\forall\var > 0$, there is an open projection $q$ such that $q\geq p$ and $\|xq\| < 1+\var$.

\proclaim{Theorem 5.7}If $p$ is a closed projection in $A^{**}$, then (4) is true for $p$ if and only if $p$ is regular in $\tilde A^{**}$.
\endproclaim

\demo{Proof}If $p$ is regular in $\tilde A^{**}$, then $\|xp\|=\|x\overline p\|$, where $\overline p$ is the closure in $\tilde A^{**}$, and we can just apply [2] for $\tilde A$.
\enddemo

If $p$ is not regular in $\tilde A^{**}$, then by 4.8, $1<\alpha (p) < \infty$.
Choose $a$ in $A_+$ such that $pap\geq p$, and let $s=\|a\|$, so that $1< s<\infty$.
Then $p(s-a)p\leq (s-1)p$.
If we apply (4) to $x=(s-1)^{-{1\over 2}}(s-a)^{1\over 2}$, we find an open projection $q$ dominating $p$ such that $q(s-a)q\leq (s-{1\over 2})q$.
Thus $qaq\geq {1\over 2}q$, so that $\alpha(q)<\infty$, and $p$ is compact by 5.4.
Since $\alpha(p) > 1$, this is impossible, and hence (4) is false for $p$.

The final result of this section is on the same subject as [8].
[8] dealt with the non-commutative analogue of open relatively compact sets (except that the correct analogue turned out to be \underbar{nearly} relatively compact projections).
Now we consider the non-commutative analogue of open sets with compact boundary.

\proclaim{Theorem 5.8}Let $A$ be a $\sigma$-unital $C^*$-algebra, $q$ an open projection in $A^{**}$, and $B=\her(q)$.
The following are equivalent:\ \ 1.\ $M(A,B)/B$ is $\sigma$-unital.

2.\ There is a closed projection $p$ such that $p\leq q$ and $\alpha(q-p)<\infty$.

3.\ $M(A,B)/B$ is unital.
\endproclaim
 
\demo{Proof}1 $\Rightarrow$ 2:\ Let $h$ be a positive element of $M(A,B)$ such that the image of $h$ is strictly positive in $M(A,B)/B$, and let $e$ be a strictly positive element of $A$.

\medskip\noindent
{\bf Claim}.\ $q(e+h)q\geq \var q$ for some $\var>0$.

The proof of the claim is similar to, and will refer to, the proof of [8, Thm.~4].
If false, we can find a sequence $(\varphi_n)$ in $P(B)$ such that $\Sigma\varphi_n (e+h)<\infty$.
We use this sequence as in [8] to construct a closed subprojection $p'$ of $q$ such that $p' Ap'\subset\cK(p' H)$ and also $p' h^{1\over 2} M(A,B) h^{1\over 2} p'\subset\cK(p' H)$.
If $C=[h^{1\over 2} M(A,B) h^{1\over 2}]^-$, then $p' Cp'\subset\cK$ and $C$ is a hereditary $C^*$-subalgebra of $m(A,B)$ whose image in $M(A,B)/B$ is everything.
Therefore $M(A,B)\subset C+B\subset C+A$.
Thus $p' M(A,B) p'\subset \cK$.
But by [7, 3.31] there is $x$ in $M(A,B)$ with $p'\leq x\leq q$.
Then $p'=p' xp'\in\cK(p' H)$, which is absurd, since $p'$ is an infinite rank projection on $H$.
Thus the claim is proved.

Now let $p=E_{[{\var\over 2},\infty) }(h)$.
Then
$$
\align
h&\leq \|h\|p+{\var\over 2} (q-p).\text{ Therefore}\\
&\var q\leq q(e+h)q,\\
&\var(q-p)\leq (q-p)(e+h)(q-p)\leq (q-p)[e+\|h\|p+{\var\over 2}(q-p)](q-p),\\
&{\var\over 2}(q-p)\leq (q-p) e(q-p).
\endalign
$$
Thus $\alpha(q-p) < \infty$.

2 $\Rightarrow$ 3:\ Choose $h$ in $M(A)$ such that $p\leq h\leq q$ ([7, 3.31]).
Then $h\in M(A,B)$.
If $x\in M(A,B)$, then $(1-h) xx^* (1-h)$ is in $M(A,\her (q-p))$.
But by [8], $M(A,\her(q-p))=\her(q-p)$.
Thus $(1-h)xx^*(1-h)$ is in $A$, which implies $(1-h)x\in A$.
Similarly, $x(1-h)\in A$.
This shows that the image of $h$ is an identity for $M(A,B)/B$.

3 $\Rightarrow$ 1 is trivial.
\enddemo

\subheading{\S 6.\ $\alpha (p_1 \vee p_2)$}

If $p_1$ and $p_2$ are closed projections with a positive angle, then $p_1\vee p_2$ is closed, by [1], but if the angle is 0, $p_1\vee p_2$ may not be closed.
The same applies to compactness, since $p$ is compact in $A^{**}$ if and only if closed in $\tilde A^{**}$.
Therefore it is natural to attempt to bound $\alpha(p_1\vee p_2)$ in terms of $\alpha(p_1),\alpha(p_2)$, and the angle between $p_1$ and $p_2$.

At the cost of some redundancy, we first prove a special case which is considerably easier than the general case and is proved differently.

\proclaim{Theorem 6.1}Assume $p_1$ and $p_2$ are projections in $A^{**}$ and $p_1 p_2=0$.\newline
(a)\ If $p_1$ and $p_2$ are closed and $A$ is $\sigma$-unital then $\alpha(p_1+p_2)=\max (\alpha (p_1),\alpha(p_2))$.\newline
(b)\ In general, $\alpha(p_1+p_2)^{-1}\geq \alpha(p_1)^{-1}+\alpha(p_2)^{-1}-1$.
\endproclaim 

\demo{Proof}(a)\ By [7, 3.31] and the continuous functional calculus, we can find $h_1,h_2$ in $M(A)_{sa}$ such that $p_j\leq h_j\leq 1$ and $h_1 h_2=0$.
If $a_j\in A_{sa}$ and $p_j a_j p_j\geq p_j$, let $b_j=h_j a_j h_j$ and $b=b_1+b_2$.
Then $\|b\|\leq\max (\|a_1\|,\|a_2\|)$, and $(p_1+p_2)b(p_1+p_2)\geq p_1+p_2$.
The result follows easily.

(b)\ This is vacuous unless $\alpha(p_1),\alpha(p_2)<\infty$.
Therefore assume this and choose $\var_1,\var_2$ such that $0<\var_j<\alpha(p_j)^{-1}$.

We use Proposition 2.8 for an approximate identity $(e_i)$, where $\|e_i\|<1$, $\forall i$.
Let $p=p_1+p_2$, and let $pe_i p$ be represented by the operator matrix $\pmatrix a_i&b_i\\ b_i^*&c_i\endpmatrix$, relative to $pH=p_1 H\oplus p_2 H$.
Since $\|a_i\|<1$, the inequality, $pe_i p\leq p$, is equivalent to $b_i^*(1-a_i)^{-1} b_i\leq 1-c_i$.
For $i$ sufficiently large $a_i\geq\var_1$ and $c_i\geq\var_2$.
Then $(1-a_i)^{-1}\geq (1-\var_1)^{-1}$, and hence $\|b\|^2\leq (1-\var_1)(1-\var_2)$.\newline
Then $\pmatrix a_i&b_i\\ b_i^*&c_i\endpmatrix\geq \pmatrix \var_1&b_i\\ b_i^*&\var_2\endpmatrix \geq\pmatrix \var_1+\var_2-1&0\\ 0&\var_1+\var_2-1\endpmatrix$,\newline
by an easy calculation, and hence $\alpha(p)^{-1}\geq\var_1+\var_2=1$.
Since $\var_j$ can be taken arbitrarily close to $\alpha(p_j)^{-1}$, the result follows.
\enddemo

\proclaim{Corollary 6.2}If $p_1$ and $p_2$ are projections in $A^{**}$ such that $p_1 p_2=0$, then:\newline
(a)\ If $\alpha(p_1),\alpha(p_2)<\infty$ and $\alpha(p_1)^{-1}+\alpha(p_2)^{-1}>1$, then $\alpha(p_1+p_2)<\infty$.\newline
(b)\ If $\alpha(p_2)=1$, then $\alpha(p_1+p_2)=\alpha(p_1)$.
\endproclaim

\medskip\noindent
{\bf Examples} 6.3.\ (a) From 6.2(b) or otherwise, we see that $p_1 p_2=0$ and $\alpha(p_1)=\alpha(p_2)=1$ imply $\alpha(p_1+p_2)=1$.
But it could be that $p_1,p_2\in RC$ and $p_1+p_2\not\in RC$.
Let $C^*(p,q)$ be the free $C^*$-algebra generated by two projections without an identity.
(This $C^*$-algebra is described in \S3 of [17].)
Let $\pi\colon C^*(p,q)\to B(H)$ be a representation which induces a one-to-one map from $C^*(p,q)$ to $B(H)/\cK$.
Let $A$ be the extension of $\cK\oplus \cK$ by $C^*(p,q)$ induced by $\pi\oplus\pi$ (cf.~4.15(a) and [10]).
Then $A^{**}$ can be identified with $B(H)\oplus B(H)\oplus C^*(p,q)^{**}$ so that any element $x$ of $C^*(p,q)$ becomes $\pi(x)\oplus\pi(x)\oplus x$.
Let $p_1=\pi(p)\oplus 0\oplus 0$ and $p_2=0\oplus\pi (q)\oplus 0$.
Since $p_1\leq p$ and $p_2\leq q$, $p_1$ and $p_2$ are relatively compact.
We claim that $(p_1+p_2)^- =\pi(p)\oplus\pi(q)\oplus 1$, a non-compact projection (since $C^*(p,q)$ is non-unital).
To prove this, we just have to show $(x+y)(p_1+p_2)=0$, $x\in C^*(p,q)$, $y\in\cK\oplus\cK$, implies $x=0$.
(It is then easy to compute $\{a\in A\colon a(p_1+p_2)=0\}$ which equals\newline
$\{a\in A\colon a(p_1+p_2)^- =0\}.)$
If $(x+y)(\pi(p)\oplus\pi(q)\oplus 0)=0$, then $\pi(x) \pi(p)$, $\pi(x)\pi(q)\in\cK$.
Therefore $\pi(x(p+q))\in\cK$, and hence $x(p+q)=0$.
Since $p+q$ is a strictly positive element of $C^*(p,q)$, $x=0$.

(b)\ We give a simple example where $\alpha(p_1)^{-1}+\alpha(p_2)^{-1}=1$ and $\alpha(p_1+p_2)=\infty$.
Let $A=c\otimes\cK$.
Choose $\theta$ in $(0,{\pi\over 2})$, and let $v_n=\cos\theta e_1+\sin\theta e_{n+1}$, $w_n=\sin\theta e_1-\cos\theta e_{n+1}$.
Define $p_1$ and $p_2$ by $(p_1)_\infty=(p_2)_\infty=0$, $(p_1)_n=v_n\times v_n$, $(p_2)_n=w_n\times w_n$.
Then $\alpha(p_1)=\cos^{-2}\theta$, $\alpha(p_2)=\sin^{-2}\theta$.
Since $(p_1+p_2)_n\geq e_{n+1}\times e_{n+1}$, $\alpha(p_1+p_2)=\infty$.

\medskip\noindent
{\bf Example 6.4}.
Before proceeding to a general result, we give a simple example to show that the hypothesis angle $(p_1,p_2) > 0$ is necessary.
Let $A=c\otimes\cK$ and $v_n=(1-n^{-1})^{1\over 2} e_1+n^{-{1\over 2}}e_{n+1}$.
Define projections $p$ and $q$ in $A^{**}$ by $p_\infty=q_\infty=e_1\times e_1$, $p_n=e_1\times e_1$, $q_n=v_n\times v_n$.
Then $p$ and $q$ are both compact.
$p\vee q$ is given by $(p\vee q)_\infty=e_1\times e_1$ and $(p\vee q)_n=e_1\times e_1+e_{n+1}\times e_{n+1}$.
Thus $p\vee q$ is closed and $\alpha(p\vee q)=\infty$.

If we consider instead $p'$ and $q'$, where $p'_n=p_n,q'_n=q_n$, and $p'_\infty=q'_\infty=0$, then we obtain disjoint, open, relatively compact, and $k$-regular projections such that $\alpha(p'\vee q')=\infty$.

For the general case, we consider two situations.

\medskip\noindent
I.\ angle$(p_1,p_2)=\theta$, $\alpha(p_j)=\sec^2\theta_j$, $0<\theta\leq {\pi\over 2}$, $0\leq\theta_j<{\pi\over 2}$, $\theta_1+\theta_2<\theta$.
Then $\alpha(p_1\vee p_2) <\infty$ and $\alpha(p_1\vee p_2)^{-1}\geq {S-\sqrt{T}\over 2\sin^2\theta}$, where
$$
\align
S&=\cos^2\theta_1+\cos^2\theta_2-2\cos^2\theta-2\cos\theta\sin\theta_1\sin\theta_2,\text{ and}\\
T&=(\cos^2\theta_1+\cos^2\theta_2)^2+4\cos^2\theta\cos^2\theta_1\cos^2\theta_2\\
&\quad -4\cos\theta(\cos^2\theta_1+\cos^2\theta_2)\sin\theta_1\sin\theta_2\\
&\quad-4(1+\cos^2\theta)(\cos^2\theta_1+\cos^2\theta_2)+8\cos\theta\sin\theta_1\sin\theta_2\\
&\quad +4\cos^2\theta+4.
\endalign
$$
If $\theta={\pi\over 2}$, this formula is the same as 6.1(b); if $\alpha(p_2)=1$, this gives $\alpha(p\vee p_2)^{-1}\geq { \alpha(p_1)^{-1}-\cos^2\theta\over\sin^2\theta}$; and if $\alpha(p_1)=\alpha(p_2)$, this gives
$$
\alpha(p_1\vee p_2)^{-1}\geq {1+\cos\theta\over\sin^2\theta} (2\alpha (p_j)^{-1} -1-\cos\theta).
$$
This estimate and the hypothesis $\theta_1+\theta_2<\theta$ are sharp, even if we add the assumption that $p_1$ and $p_2$ are disjoint, open, and $k$-regular, $\forall k$, or if we add the assumption that $p_1$ are $p_2$ are closed.

\medskip\noindent
II.\ $p_1$ and $p_2$ are closed \underbar{and} $p_1\wedge p_2=0$, angle $(p_1,p_2)=\theta$, $\alpha(p_j)=\sec^2\theta_j$, $0<\theta\leq {\pi\over 2}$, $0\leq \theta_j< {\pi\over 2}$.
Then $\alpha(p_1\vee p_2)<\infty$ and:
\item{(a)}If $\cos\theta\leq {\sin\theta_1\sin\theta_2\over 1+\cos\theta_1\cos\theta_2}$, then $\alpha(p_1\vee p_2)^{-1}\geq {S-\sqrt{T}\over 2\sin^2\theta}$, where
$$
\aligned
S&=\cos^2\theta_1 +\cos^2\theta_2+2\cos^2\theta\cos\theta_1\cos\theta_2,\text{ and}\\
T&=(\cos^2 \theta_1+\cos^2\theta_2)^2+4\cos^2\theta_1\cos^2 \theta_2 (\cos^2\theta-\sin^2\theta)\\
&\quad +4\cos\theta_1\cos\theta_2 (\cos^2\theta_1+\cos^2\theta_2)\cos^2\theta
\endaligned
$$
\item{(b)}If $\cos\theta\geq {\sin\theta_1\sin\theta_2\over 1+\cos\theta_1\cos\theta_2}$, then
$$
\alpha(p_1\vee p_2)^{-1}\geq {\cos^2\theta_1\cos^2\theta_2\sin^2\theta\over \cos^2\theta_1+\cos^2\theta_2-\cos^2\theta_1\cos^2\theta_2(1+\cos^2\theta)+2\cos\theta\cos\theta_1\cos\theta_2\sin\theta_1\sin\theta_2}
$$
If $\theta={\pi\over 2}$, this gives $\alpha(p_1\vee p_2)=\max(\alpha(p_1),\alpha(p_2))$; if $\alpha(p_2)=1$, this gives $\alpha(p_1\vee p_2)\leq {\alpha(p_1)-\cos^2\theta\over\sin^2\theta }$; and if $\alpha(p_1)=\alpha(p_2)$, this gives
$$
\alpha(p_1\vee p_2)\leq\cases {1+\cos\theta\over 1-\cos\theta}\ \alpha(p_j),&\cos\theta\leq {\alpha(p_j)-1\over\alpha(p_j)+1}\\
{2-\alpha(p_j)^{-1} (1+\cos\theta)\over 1-\cos\theta}\ \alpha(p_j),&\cos\theta\geq {\alpha(p_j)-1\over\alpha(p_j)+1}.\endcases
$$
This estimate is sharp.

There are some preliminaries before the proof of the positive results.
First, the angle between $p_1$ and $p_2$ is the same as the angle between $p_1$ and $p_2-p_1\wedge p_2$.
Thus in both cases, we may assume $p_1\wedge p_2=0$.
Then if $\varphi\in F(p_1\vee p_2)$, we can write $\varphi=\sum\limits^2_{j,k=1}\ \varphi^{jk}$, where $\varphi^{jk}\in \{f\in A^*\colon f(\cdot)=f(p_j\cdot p_k)\}$.
We do this by considering $A^{**}$ as a subalgebra of $B(H)$ via the universal representation of $A$.
Then $\varphi=(\cdot v,v), v\in (p_1\vee p_2)H$.
Since angle $(p_1,p_2)>0$, $(p_1\vee p_2)H=p_1 H+p_2 H$, and $v=v_1+v_2$, $v_j\in p_j H$.
Then $\varphi^{11}=(\cdot v_1,v_1)$, $\varphi^{12}=(\cdot v_2,v_1)$, etc.
Note that $\|\varphi\|=\varphi(1)=\varphi^{11}(1)+2R e\varphi^{12} (1)+\varphi^{22}(1)=\|\varphi^{11}\|+2Re\varphi^{12}(1)+\|\varphi^{22}\|$, and $|\varphi^{12}(1)|\leq \|\varphi^{11}\|^{1\over 2}\|\varphi^{22}\|^{1\over 2}\cos\theta$.
It is important to know that the $\varphi^{jk}$'s are uniquely determined by $\varphi=\Sigma\varphi^{jk}$.
To see this, note that \newline
$\varphi^{11}+\varphi^{21}\in L(p_1)=\{f\in A^*\colon f(\cdot)=f(\cdot p_1)\}$, and $\varphi^{12}+\varphi^{22}\in L(p_2)$.
It is easy to see that $p_1\wedge p_2=0$ implies $L(p_1)\cap L(p_2)=0$.
The reader can easily complete the proof that the four vector spaces in the decomposition are linearly independent.
Finally, we will use a slightly different notation in the actual proof.
Write $v_1=su_1$, $v_2=tu_2$, where $s,t\geq 0$ and $\|u_1\|=\|u_2\|=1$.
Then let $\psi^{11}=(\cdot u_1,u_1)$, $\psi^{12}=(\cdot u_2,u_1)$, etc., so that $\varphi=s^2 \psi^{11}+2 stRe\psi^{12}+t^2\psi^{22}$.
Note that the hypothesis angle $(p_1,p_2)=\theta$ implies $|\psi^{12}(1)|\leq\cos\theta$, and this implies $s^2+t^2\leq (1-\cos\theta)^{-1}$.

One more remark may be helpful.
In the proof below we first show that $\alpha(p_1\vee p_2)^{-1}$ is at least the solution to a certain minimum problem for a function of several real variables.
We then sketch the solution of this minimum problem.
In the examples where we show our bounds are sharp, we use the minimum problem itself rather than the explicit formula.
Thus the reader may not wish to verify that our solution of the minimum problem is correct.

\proclaim{Theorem 6.5}The upper bounds given for $\alpha(p_1\vee p_2)$ in I and II above are valid under the hypotheses stated.
\endproclaim

\demo{Proof}We use Theorem 2.9.
Thus let $\varphi_i\in F(p_1\vee p_2)\cap S(A)$ and assume $\varphi_i\overset{w^*}\to\longrightarrow\varphi$.
Of course, in case I, $\varphi$ may not be in $F(p_1\vee p_2)$.
Using the notation above and passing to a subnet, we may assume $s_i\to s,\ t_i\to t,\ Re\psi_i^{12} (1)\to x$, and $\psi_i^{jk}\overset{w^*}\to\rightarrow \psi^{jk}$.
Let $y=Re\psi^{12} (1)$.
($y$ need not equal $x$, since $1\not\in A$.)
Also let $\delta_j=\psi^{jj}(1)$.
Clearly, $s^2+2stx+t^2=1$, $|x|\leq\cos\theta$, and $\delta_j\geq\cos^2\theta_j$.
Also, by the lower semicontinuity of norm\newline
$\|(s')^2 \psi^{11}+2s' t' Re \psi^{12} +(t')^2\psi^{22}\|\leq\liminf\|(s')^2 \psi_i^{11}+2s' t' Re\psi_i^{12}+(t')^2 \psi_i^{22}\|$; i.e., \newline
$\delta_1(s')^2+2ys' t' +\delta_2 (t')^2\leq (s')^2 +2xs' t'+(t')^2$, $\forall s'$, $t'\in\bR$.\newline
Thus $\pmatrix 1-\delta_1 & x-y\\ x-y & 1-\delta_2\endpmatrix\geq 0$, and $|x-y|\leq (1-\delta_1)^{1\over 2} (1-\delta_2)^{1\over 2}\leq\sin\theta_1\sin\theta_2$.
Since $\varphi=s^2\psi^{11}+2st Re\psi^{12}+t^2 \psi^{22}$, we find that $\|\varphi\|$ is at least the minimum of\newline
$\cos^2\theta_1 s^2+2yst+\cos^2\theta_2 t^2$ subject to $s^2+2xst+t^2=1$, $|x|\leq\cos\theta$, $|x-y|\leq\sin\theta_1\sin\theta_2$, and $s,t\geq 0$.

For case I, we compute this minimum and show it is the formula given.
Note that $\theta_1+\theta_2 < \theta$ implies $\cos\theta_1 \cos\theta_2-\sin\theta_1\sin\theta_2>\cos\theta$.
Thus $y\geq -\cos\theta-\sin\theta_1\sin\theta_2>-\cos\theta_1\cos\theta_2$.
Thus the minimum \underbar{is} positive.
One can see without computation that at the minimum $x=-\cos\theta$ and $y=-\cos\theta-\sin\theta_1\sin\theta_2$.
(It is obvious that $y=x-\sin\theta_1\sin\theta_2$.
To see that $x=-\cos\theta$, note that if $x$ and $y$ are decreased by the same amount (for fixed $s,t$), both quadratics change by the same amount, and thus the smaller quadratic changes by the larger percentage.)
Once $x$ and $y$ are known, it is a matter of routine calculus (Lagrange multipliers) to calculate the minimum; and this will be left to the reader.

In case II, $\varphi\in F(p_1\vee p_2)$ and $\psi^{jk}(\cdot)=\psi^{jk}(p_j\cdot p_k)$.
Thus, using the same notation as for case I, we find that $\|\varphi\|=\delta_1 s^2+2yst+\delta_2 t^2$ and $y\leq\delta_1^{1\over 2}\delta_2^{1\over 2}\cos\theta$.
Thus now $\|\varphi\|$ is at least the minimum of $\delta_1 s^2+2yst+\delta_2 t^2$ subject to $s^2+2xst+t^2=1$, $|x|\leq\cos\theta$, $|y|\leq \delta_1^{1\over 2}\delta_2^{1\over 2}\cos\theta$, $|x-y|\leq (1-\delta_1)^{1\over 2} (1-\delta_2)^{1\over 2}$, $\cos^2\theta_j\leq\delta_j\leq 1$, and $s,t\geq 0$.
(Unlike case I, it is not yet obvious that $\delta_j=\cos^2\theta_j$ at the minimum.)

We can see by reasoning similar to that of case I that for fixed $\delta_1,\delta_2$, the minimum occurs at $y=-\delta_1^{1\over 2}\delta_2^{1\over 2}\cos\theta$ and 
$$
x=\cases y+(1-\delta_1)^{1\over 2}(1-\delta_2)^{1\over 2},&(1-\delta_1)^{1\over 2}(1-\delta_2^{1\over 2})\leq\cos\theta(1+\delta_1^{1\over 2}\delta_2^{1\over 2})\\
\cos\theta,&(1-\delta_1)^{1\over 2} (1-\delta_2)^{1\over 2}\geq\cos\theta (1+\delta_1^{1\over 2}\delta_2^{1\over 2}).\endcases
$$
We then substitute these values of $x$ and $y$ and prove that the minimum in $(s,t)$ is a monotone increasing function of $\delta_1$ and $\delta_2$ (so that the minimum occurs for the smallest values of $\delta_1,\delta_2$).
The easiest way to see the monotonicity is to perform a change of variable:\newline
Replace $s$ by $\delta_1^{-{1\over 2}}s$ and $t$ by $\delta_2^{-{1\over 2}}t$.
The rest of the calculation is left to the reader.

\medskip\noindent
{\bf Remark}.
It follows from the formulas in both I and II that $\alpha(p_1)=\alpha(p_2)=1$ implies $\alpha(p_1\vee p_2)=1$, but this is nothing new.
It follows from $p_1\vee p_2\leq K(\theta) (p_1+p_2)$, where, of course, $K(\theta)\to\infty$ as $\theta\to 0$.
\enddemo

\proclaim{Corollary 6.6}If $p_1$ and $p_2$ are closed projections such that $p_1\wedge p_2=0$ and at least one of $\alpha(p_1),\alpha(p_2)$ is finite, then angle $(p_1,p_2)>0$.
Thus if both of $\alpha(p_1),\alpha(p_2)$ are finite, then $\alpha(p_1\vee p_2)$ is finite.
In particular, if $p_1$ and $p_2$ are compact and $p_1\wedge p_2=0$, then $p_1\vee p_2$ is compact.
\endproclaim

\demo{Proof}Assume angle $(p_1,p_2)=0$.
Then there are unit vectors $v_n$ in $p_1 H$ and $w_n$ in $p_2 H$ such that $\|v_n-w_n\|\to 0$.
Thus there are states $\varphi_n$ in $F(p_1)$ and $\psi_n$ in $F(p_2)$ such that $\|\varphi_n-\psi_n\|\to 0$.
Assume $\alpha(p_1) <\infty$, and let $\varphi$ be a weak$^*$ cluster point of $(\varphi_n)$.
Then $\varphi\neq 0$ and $\varphi$ is also a weak$^*$ cluster point of $(\psi_n)$.
Thus $F(p_1)\cap F(p_2)\neq \{0\}$, a contradiction.
\enddemo

\medskip\noindent
{\bf Examples 6.7}.
We show the sharpness claimed in I and II.
Let $A=c\otimes\cK$.

I.\ Let $\theta,\theta_1$ and $\theta_2$ be as above, except that now we allow the possibility that $\theta_1+\theta_2=\theta$.
Choose $x,y,s,t$ in $\bR$ such that $s^2+2xst+t^2=1$, $|x|\leq\cos\theta$, $|x-y|\leq\sin\theta_1\sin\theta_2$, $s,t\geq 0$, and $\cos^2\theta_1 s^2+2yst+\cos^2 \theta_2 t^2$ is minimized subject to the above.
Of course, we know that $x=-\cos\theta$ and $y=-\cos\theta-\sin\theta_1\sin\theta_2$, and we could calculate $s,t$.
If $\theta_1+\theta_2=\theta$, it is easily seen that this minimum value is 0, and hence the example in this case will have $\alpha(p_1\vee p_2)=\infty (0\in\overline S(p_1\vee p_2))$.

Choose vectors $u^1,u^2$ in $H$ such that $\|u^j\|=\cos\theta_j$ and $(u^1,u^2)=y$.
The proof of 6.5 showed that $|y|\leq\cos\theta_1\cos\theta_2$ (actually equality), and therefore this is possible.
For each $n$ choose vectors $w_n^1,w_n^2$ in $H$ such that $\|w_n^j\|=\sin\theta_j$, $(w_n^j,u^k)=0$, $(w_n^1,w_n^2)=x-y$, and $w_n^j\overset w\to\rightarrow 0$ as $n\to\infty$.
This is clearly possible.
Let $u_n^j=u^j+w_n^j$, and let $v_n=su_n^1+tu_n^2$.
Then $\|u_n^j\|=1$, $(u_n^1,u_n^2)=x$ and $u_j^n\overset w\to\rightarrow u^j$ as $n\to\infty$.
It follows that $\|v_n\|=1$.
Let $r$ be the projection in $B(H)$ whose range is span$(u^1,u^2)$.

Define closed projections $p^1,p^2$ in $A^{**}$ by $p_\infty^1=p_\infty^2=r$, $p_n^j=u_n^j\times u_n^j$.
Define $a$ in $A_{sa}$ by $a_n=a_\infty=r$.
Then $\|a\|=1$ and $p^j ap^j\geq\cos^2\theta_j p^j$.
Therefore $\alpha(p^j)^{-1}\geq\cos^2\theta_j$.
Since angle$(p^1,p^2)=$ angle$(p^1-p^1\wedge p^2,p^2-p^1\wedge p^2)$, and since $|(u_n^1,u_n^2)|\leq\cos\theta$, angle$(p^1,p^2)\geq\theta$.
Let $\varphi_n$ in $S(A)$ be defined by $\varphi_n(a)=(a_n v_n,v_n)$.
Clearly $\varphi_n\in F(p^1\vee p^2)$ and $\varphi_n \overset{w^*}\to\rightarrow\varphi$, where $\varphi(a)=(a_\infty v,v)$, $v=su^1+tu^2$.
Therefore $\|\varphi\|=\|v\|^2=\cos^2 \theta_1 s^2+2yst+\cos^2\theta_2 t^2$.
This shows that $\alpha(p^1\vee p^2)$ is at least the value specified in I.
(Of course, by 6.5, the inequalities for $\alpha(p^j),\alpha(p^1\vee p^2)$, and angle$(p^1,p^2)$ are actually equalities.)

Now we show how to modify the above to obtain disjoint, open, $k$-regular projections $q^1,q^2$.
Let $q_\infty^j=0$,
$$
\align
q_n^1&=\cases p_m^1,&n=3m\\
{u^1\over \|u^1\|}\times{u^1\over \|u^1\|},&n=3m+1\\
0,&n=3m+2,\endcases\\
\noalign{\medskip}
\text{and }q_n^2&=\cases p_m^2,&n=3m\\
0,&n=3m+1\\
{u^2\over \|u^2\|}\times {u^2\over \|u^2\|},&n=3m+2.\endcases
\endalign
$$
The reader can easily verify that $q^1,q^2$ have the required properties.
(The closures of $q^j$ have $(\overline q^j)_\infty={u^j\over \|u^j\|}\times {u^j\over \|u^j\|}$.)

II.\ The construction is very similar.
Of course now we have hardly any restrictions on $\theta,\theta_1,\theta_2$.
Choose $x,y,s,t,\delta_1,\delta_2$ in $\bR$ such that $s^2+2xst+t^2=1$, $|x|\leq\cos\theta$, $|y|\leq \delta_1^{1\over 2}\delta_2^{1\over 2}\cos\theta$, $|x-y|\leq (1-\delta_1)^{1\over 2} (1-\delta_2)^{1\over 2}$, $\cos^2\theta_j\leq\delta_j\leq 1$, $s,t\geq 0$ and $\delta_1 s^2+2yst+\delta_2 t^2$ is minimized subject to the above.
Of course we know that $\delta_j=\cos^2 \theta_j$, $y=-\delta_1^{1\over 2}\delta_2^{1\over 2}\cos\theta$, the formula for $x$ was given in the proof of 6.5, and we could calculate $s,t$.

The definitions of $u^j,w_n^j,u_n^j$ and $v_n$ are by the same formulas used in part I except that now $\|u^j\|=\delta_j^{1\over2}$ and $\|w_n^j\|=(1-\delta_j)^{1\over 2}$.

Disjoint closed projections $p^1,p^2$ in $A^{**}$ are defined by $p_\infty^j={u^j\over \|u^j}\|\times {u^j\over \|u^j\|}$, and $p_n^j=u_n^j\times u_n^j$.
It is easy to see that angle$(p^1,p^2)\geq\theta$.
Everything else is the same as in part I.

\subheading{\S 7.\ Attainment of Extreme Values}

If $p$ is a projection in $A^{**}$, we say that $\alpha(p)$ is attained if $\alpha(p)<\infty$ and there is $a$ in $A_{sa}$ such that $\|a\|=\alpha(p)$ and $pap\geq p$.
We say dist$(p,RC)$ is attained if dist$(p,RC)<1$ and there is $q$ in $RC$ such that $\|p-q\|=$ dist$(p,RC)$.
We define attainment similarly for dist$(P,ORC)$ and dist$(p,CRC)$.

\proclaim{Proposition 7.1}Let $p$ be a projection in $A^{**}$.

(a)\ If $\alpha(p)=1$, then $\alpha(p)$ is attained if and only if dist$(p,RC)$ is attained if and only if $p\in RC$.

(b)\ If $p$ is closed, then $\alpha(p)$ is attained if and only if dist$(p,RC)$ is attained if and only if dist$(p,CRC)$ is attained.

(c)\ If $p$ is open, then dist$(p,RC)$ is attained if and only if dist$(p,ORC)$ is attained.

(d)\ In general, if dist$(p,RC)$ is attained, then $\alpha(p)$ is attained.
\endproclaim

\demo{Proof}(d)\ By the proof of Theorem 2.2, if $q$ is projection in $RC$ such that $\|p-q\|=$ dist$(p,RC)$, then $pqp\geq\alpha(p)^{-1}p$.
Let $a_1$ be in $A_{sa}$ such that $q\leq a_1\leq 1$, and $a=\alpha(p)a_1$.
Then $\|a\|=\alpha(p)$ and $pap\geq p$.

(a)\ The second equivalence is obvious, since dist$(p,RC)=0$.
In view of (d), we need just assume $\alpha(p)$ is attained and prove $p$ is in $RC$.
If $a$ is in $A_{sa}$, $\|a\|\leq 1$, and $pap\geq p$, the proof of Theorem 2.1 shows that $ap=pa$.
Therefore $p\leq a_+\leq 1$, and $p$ is in $RC$.

(b)\ Since the two distances are the same, dist$(p,CRC)$ attained implies dist$(p,RC)$ attained.
Theorem 2.6(b) and 7.1(d) complete the proof.

(c)\ Again the two distances are the same.
Thus assume dist$(p,RC)$ is attained.
Choose $q$ in $RC$ such that $pqp\geq\alpha(p)^{-1}p$, as in the proof of (d) and choose $a$ in $A_{sa}$ such that $q\leq a\leq 1$.
Choose a continuous function $f\colon\bR\to [0,1]$ such that $f(1)=1$ and $f=0$ on $(-\infty,{1\over 2}]$, and let $b=f(a)$.
Then the range projection of $b$ is in $RC$, and $pbp\geq\alpha(p)^{-1}p$ since $b\geq q$.
By the proof of Theorem 2.6(a), the range projection of $b^{1\over 2}p$ is in $ORC$ and it attains dist$(p,ORC)$.
\enddemo

\noindent
{\bf Examples 7.2}.\ (a) An open projection $p$ such that $1<\alpha (p) < \infty$ and $\alpha(p)$ is not attained.
Let $A=\{a\in c\otimes {\cK}\colon a_\infty$ is diagonal$\}$.
Then $A^{**}=\{h\in ( c\otimes{\cK})^{**}\colon h_\infty$ is diagonal$\}$, where in both cases diagonality is with respect to our usual fixed orthonormal bases of $H$.
Let $v_0$ be a unit vector in $H$ with all coordinates non-zero.
Let $\{f_1,f_2,\ldots\}$ be an orthonormal basis for $\{v_0\}^\perp$, and let $v_n=2^{-{1\over 2}}v_0+2^{-{1\over 2}}f_n$.
Define $p$ by $p_\infty=0$ and $p_n=v_n\times v_n$.
Define $x'$ in $A_{sa}$ by $x'_n=x'_\infty=\Sigma_1^k e_i\times e_i$.
Then $p_n x'_n p_n=\var_n p_n$ where $\lim\limits_{n\to\infty}\var_n={1\over 2}\sum_1^k|(v_0,e_i)|^2$.
For any $\delta>0$, we can modify $x'_n$ for finitely many values of $n$ to obtain $x$ in $A_{sa}$ such that $\|x\|=1$ and $pxp\geq (2^{-1}\sum_1^k|(v_0,e_i)|^2-\delta)p$.
Since $k$ can be arbitrarily large and $\delta$ arbitrarily small, $\alpha(p)^{-1}\geq 2^{-1}$.

We claim there is no $x$ in $A_{sa}$ such that $\|x\|\leq 1$ and $pxp\geq 2^{-1}p$.
If such $x$ existed, $p_n x_\infty p_n\geq (2^{-1} -\delta_n)p_n$, where $\delta_n\to 0$ as $n\to\infty$.
Since $v_n\overset w\to\rightarrow 2^{-{1\over 2}}v_0$, this implies $2^{-1}(x_\infty v_0,v_0)\geq 2^{-1}$.
This is impossible for $x_\infty$ compact and diagonal.
It follows from the proved claim that $\alpha(p)=2$.

(b)\ An open projection $p$ such that $\alpha(p)$ is attained and dist$(p,RC)$ is not attained.
Let $A=c\otimes\cK$.
Let $a_0=$ Diag$(1,d_2,d_3,\ldots)$, where $0<d_n<{1\over 2}$, $d_n\to 0$, and $\sum_2^\infty d_n=\infty$.
For $n\geq 2$, let $v_n=\left({ 2^{-1}-d_n\over 1-d_n}\right)^{1\over 2} e_1+\left( { 2^{-1}\over 1-d_n }\right)^{1\over 2} e_n$.
Then $(a_0 v_n,v_n)={1\over 2}$ and $\|v_n\|=1$.
Define $p$ so that $p_\infty=0$ and the sequence $(p_n)$ includes each $v_k\times v_k$ infinitely often.
Define $a$ in $A_{sa}$ by $a_\infty=a_n=a_0$.
Then $\|a\|=1$ and $pap\geq {1\over 2} p$.
We claim there is no $q$ in $RC$ such that $pqp\geq {1\over 2} p$.
If there is, there is $c$ in $A_{sa}$ such that $q\leq c\leq 1$.
Then $q\leq E_{\{1\}}(c)$.
It is easy to see from the holomorphic functional calculus that $E_{\{1\}}(c)$ is majorized by a projection in $A$.
Changing notation, we assume $q\in A$.
For any $k$, we can choose $n_1<n_2<\ldots$, such that $p_{n_i}=v_k\times v_k$.
Then $(q_{n_i} v_k,v_k)\geq {1\over 2}$ and $q_n\to q_\infty$ imply $(q_\infty v_k,v_k)\geq {1\over 2}$.
Let $r=q_\infty\vee (e_1\times e_1)$.
Thus $r$ is a finite rank projection, $r=e_1\times e_1+r'$ for a projection $r'$, and $(rv_k,v_k)\geq {1\over 2}$, $\forall k$.
Now $(rv_k,v_k)={2^{-1}-d_k\over 1-d_k}+{2^{-1}\over 1-d_k} (r'e_k,e_k)\geq{1\over 2}$ implies $(r' e_k,e_k)\geq d_k$.
But $r'$ finite rank implies $\sum_1^\infty (r' e_k,e_k)<\infty$ in contradiction to the choice of $(d_k)$.
Again we can conclude \underbar{a postiori}, that $\alpha(p)=2$.

(c)\ A closed projection $p$ such that $1<\alpha(p)<\infty$ and $\alpha(p)$ is not attained.
Let $A_0$ be the $C^*$-algebra called $A$ in (a), and let $A_1=A_0\otimes M_2$.
Let $e'$ be the projection in $M(A_1)$ given by $\pmatrix {1\over 2}&{1\over 2}\\ {1\over 2}&{1\over 2}\endpmatrix$; and as in 3.8 and 3.9, let $A$ be the extension of $A_1$ by $\bC$ induced by $e'$, and let $e$ be the corresponding projection in $A$.
Let $v'_0$ be a unit vector in $H$ with all coordinates non-zero, and let $v_0=v'_0\oplus 0$ in $H\oplus H$.
Choose an orthonormal sequence $f_1,f_2,\ldots$, in $H\oplus H$ such that $(f_n,v_0)=0$ and $e' f_n=0$, $\forall n$.
Let $v_n=2^{-{1\over 2}} v_0+2^{-{1\over 2}}f_n$, and define a closed projection $p'$ in $A_1^{**}$ by $p'_n=v_n\times v_n$ and $p'_\infty=\pmatrix 1&0\\ 0&0\endpmatrix$.
Identify $A^{**}$ with $A_1^{**}\oplus\bC$ and define $p=p'\oplus 1$, so that $p$ is a closed projection in $A^{**}$.
For each $m$ let $Q_m=\sum_1^m e_k\times e_k$, a diagonal projection in $\cK$, and define $a_m$ in $A_1$ by $(a_m)_n=(a_m)_\infty=\pmatrix {1\over 2}Q_m&-{1\over 2}Q_m\\ -{1\over 2}Q_m&{1\over 2}QA_m\endpmatrix$.
Then define $b_m=e+a_m$, an element of $A$.
Note that for $1\leq n\leq\infty$, the $n$'th component of $e'+a_m$ is $\pmatrix Q_m+{1\over 2}(1-Q_m)&{1\over 2}(1-Q_m)\\ {1\over 2}(1-Q_m)&Q_m+{1\over 2}(1-Q_m)\endpmatrix$, a projection.

We claim that $pb_m p\geq\var_m p$, where $\var_m\to {1\over 2}$ as $m\to\infty$.
Thus $\alpha(p)^{-1}\geq {1\over 2}$.
To prove this, it is enough to consider the $A_1^{**}$-components.
It is obvious that $p'_\infty(e'+a_m)_\infty p'_\infty\geq {1\over 2} p'_\infty$.
For $n$ finite $p'_n(e'+a_m)_n p'_n=\var_{nm} p'_n$, where\newline
$\var_{nm}={1\over 2}((e'+a_m)_n v_0,v_0)+Re((e'+a_m)_n v_0,f_n)+{1\over 2}((e'+a_m)_n f_n,f_n)$.
For all $n$, $\|(e'+a_m)_n v_0-v_0\|=2^{-{1\over 2}}\|(1-Q_m) v_0\|=\delta_m$, where $\delta_m\to 0$.
Thus\newline
$\var_{nm}\geq {1\over 2}(1-\delta_m)-(\delta_m+|(v_0,f_n)|)\geq {1\over 2}-{3\over 2}\delta_m\ \forall n$.

Now let $\varphi_n$ be the state given by $\varphi_n(x)=(x_n v_n,v_n)$.
Since $v_n\overset w\to\longrightarrow 2^{-{1\over 2}} v_0,\ \varphi_n$ converges on $A_1$ to ${1\over 2}\varphi$, where $\varphi$ is the state on $A_1$ defined by $\varphi(x)=(x_\infty v_0,v_0)$.
Also, since $e' f_n=0$, $\varphi_n(e)=(e' v_n,v_n)={1\over 2}(e' v_0,v_0)={1\over 4}$.
Thus $\varphi_n$ converges in the weak$^*$ topology of $A^*$ to ${1\over 2}\varphi\oplus 0$, where $A^*$ is identified, as usual, with $A_1^*\oplus\bC$.
This shows not only that $\alpha(p)^{-1}\leq {1\over 2}$ but also that if $a\in A_{sa}$ and $pap\geq {1\over 2}p$, then $\varphi(a)\geq 1$.
If $\alpha(p)$ were attained, then there would be such $a$ with $\|a\|=1$.
If the $A_1^{**}$-component of $a$ has $\infty$-component $\pmatrix r&*\\ *&*\endpmatrix$, then $r\leq 1$ and $r={1\over 2}\lambda+K$, where $\lambda\leq 1(a=\lambda e+a_1,a_1\in A_1)$ and $K$ is a diagonal compact operator.
Since $\varphi(a)=1$, $(rv'_0,v'_0)=1$; and this is impossible, since $r$ is diagonal, all components of $v'_0$ are non-zero, and $r\neq 1$.
Thus $\alpha(p)$ is not attained.

\subheading{\S 8.\ Majorization, $\alpha(p)$, and Semicontinuity}

If $h\in A_{sa}^{**}$ and $h\leq a$ for some $a$ in $A_{sa}$, then we expect that some of the spectral projections of $h$ will be nearly relatively compact.
This is trivially proved and has a couple of complements, one related to semicontinuity.

\proclaim{Proposition 8.1}If $h\leq a$, where $h\in A_{sa}^{**}$ and $a\in A_{sa}$, then $\alpha(E_{[\var,\infty)}(h))\leq {\|a\|\over\var}$, $0<\var\leq \|h_+\|$.
\endproclaim

\demo{Proof}If $p=E_{[\var,\infty)}(h)$, then $pap\geq php\geq\var p$.
\enddemo

\proclaim{Corollary 8.2}If, in addition, $h\geq 0$, then $\alpha(E_{[\var,\infty)}(h))=1$.
\endproclaim

\demo{Proof}Now we have $p\leq\var^{-1}/h\leq\var^{-1}a$.
\enddemo

\proclaim{Corollary 8.3}If $h\in (A_{sa})_m^-$ and $h_+\neq 0$, then $\alpha(E_{[\var,\infty) }(h))\leq {\|h_+\|\over\var}$\newline
for $0<\var\leq \|h_+\|$.
Also $E_{ \{\|h_+\| \}}(h)$ is compact.
Similarly, if $h\in\overline{A_{sa}^m}$ and $h_-\neq 0$, then $\alpha(E_{(-\infty , -\var] }(h))\leq {\|h_-\|\over\var}$ for $0<\var\leq \|h_-\|$ and
$E_{ \{-\|h_-\|\}}(h)$ is compact.
\endproclaim

\demo{Proof}Since $h$ is strongly usc, [7, 3.16] implies there is $a$ in $A_{sa}$ such that $h\leq a\leq \|h_+\|$.
Thus the inequality follows from 8.1.
Since $\|h_+\|-h$ is positive and strongly lsc, [7, 2.44(a)] implies that its range projection is open; in other words, $E_{ \{h_+\|\}}(h)$ is closed.
Also $\alpha(E_{\{\|h_+\|\}}(h))=1$, by the case $\var=\|h_+\|$, and hence this projection is compact.
The second part follows from the first applied to $-h$.
\enddemo

\medskip\noindent
{\bf Remarks}. 
2.4 and the compactness assertion of 8.3 are closely analogous to 2.44 of [7].
(Before going on, we should remind the reader that in [7] we disclaimed originality for 2.44 and much of the rest of \S2.D.)
If $p$ is a projection in $A^{**}$, then $p$ is open if and only if lsc (in any sense) ([5]), closed if and only if weakly or middle usc ([5]), and compact if and only if strongly usc ([7, 2.47]).
With the help of 2.4 and 8.3, we can now state a symmetrical result containing this last:\ If $h\in A_{sa}^{**}$ and $\sigma(h)$ has at most two elements, then $h$ is weakly or middle lsc if and only if $q$-lsc, and $h$ is strongly lsc if and only if strongly $q$-lsc.
Of course, the same is true for usc.
We also mention that there are at least two other ways of proving the compactness assertion in 8.3.
The other proofs would not mention $\alpha(p)$.

\proclaim{Corollary 8.4}Assume $h\in A_+^{**}$, $h$ is strongly usc, and $h$ is $q$-usc.
Then $h$ is strongly $q$-usc.
\endproclaim

\demo{Proof}Let $p=E_{[\var,\infty)}(h)$ for $\var$ in $(0,\|h\|]$.
By [7.3.22], $h\leq a$ for some $a$ in $A_{sa}$.
Then 8.2 implies $\alpha(p)=1$.
By the definition of $q$-usc, $p$ is closed.
Therefore $p$ is compact.
Then by definition, $h$ is strongly $q$-usc.
\enddemo

\medskip\noindent
{\bf Example 8.5}.
This example will show that the positivity assumption in 8.4 is necessary and that the estimate for $\alpha(E_{[\var,\infty) }(h))$ in 8.3 is sharp.
Also $\sigma(h)$ has only three elements.
Choose $\lambda_1,\lambda_2$ in $\bR$ such that $\lambda_1 > 2\lambda_2 > 0$.
Let $A=c\otimes\cK$ and $v_n=2^{-{1\over 2}} e_1+2^{-{1\over 2}} e_{n+1}$.
Define $p$ in $A^{**}$ by $p_\infty=e_1\times e_1$ and $p_n=v_n\times v_n$.
Also define $p_0$ by $(p_0)_\infty=e_1\times e_1$ and $(p_0)_n=0$.
Then $p_0$ is a compact projection, $p$ is a closed projection, $\alpha(p)=2$, and $p_0\leq p$.
Let $h$ in $A_{sa}^{**}$ be determined by $\sigma(h)=\{\lambda_1,\lambda_2,\lambda_3\}$, $E_{\{\lambda_1\}} (h)=p_0$, and $E_{\{\lambda_1,\lambda_2\}} (h)=p$, where $\lambda_3$ is a negative number to be determined.
Clearly, $h$ is $q$-usc, and since $p$ is not compact, $h$ is not strongly $q$-usc.
[7, 5.13 and Remark (i)] gives the following criterion for determining that $h$ is strongly usc:\ Choose a sequence $(k_m)$ in $\cK$ such that $k_m\searrow h_\infty$.
Then we require that for each $m$ and each $\var>0$, there is $N$ such that $k_m\geq h_n-\var$ for $n\geq N$.
We can take $k_m=\lambda_1 e_1\times e_1+\lambda_3\sum_2^m e_k\times e_k$.
Then if $n\geq m$, we need look only at span$(e_1,e_{n+1})$ to check the inequality $k_m\geq h_n-\var$.
It is sufficient that \newline
$\pmatrix \lambda_1&0\\ 0&0\endpmatrix\geq \lambda_2 \pmatrix {1\over 2}&{1\over 2}\\ {1\over 2}&{1\over 2}\endpmatrix +\lambda_3 \pmatrix {1\over 2}&-{1\over 2}\\ -{1\over 2}&{1\over 2}\endpmatrix$.
The reader can easily check that this is true for $|\lambda_3|$ sufficiently large.

If $\var=\lambda_2$, the inequality of 8.3 states that $\alpha(p)\leq {\lambda_1\over\lambda_2}$.
Since ${\lambda_1\over\lambda_2}$ can be close to 2, the estimate in 8.3 cannot be improved.

By slightly modifying this example, we can show that the inequality in 8.3 is not valid under the weaker hypothesis of 8.1.
Let $h'$ in $A_{sa}^{**}$ be determined by $\sigma(h')=\{\lambda_2,\lambda_3\}$ and $E_{\{\lambda_2\}}(h')=p$.
Then $h'$ is $q$-usc, and $h'$ satisfies all the hypothesis of 8.3 except that it is not strongly usc.

\subheading{\S9.\ Concluding Remarks}

\item{1.}The reader has probably noticed that in many of our examples $p$ is abelian, in the (usual) sense that the $W^*$-algebra $pA^{**}p$ is abelian.
In a few examples $\overline p$ is also abelian.
We have not systematically tried to determine which phenomena can be exhibited with abelian projections.
We merely were making a reasonable effort to keep our examples simple.
It might be interesting to know the consequences of the hypothesis $p$ is abelian or the hypothesis $\overline p$ is abelian.

\item{2.}The idea of looking at dist$(p,CRC)$ for general projections $p$ was an afterthought.
Of course, if $p$ is closed, dist$(p,CRC)=$ dist$(p,RC)$, and if $p$ is open but not closed dist$(p,CRC)=1$ by 5.1.
For $p$ neither open nor closed, all we have done is to look at the most obvious example.

Let $A=c\otimes\cK$ and $v_n=2^{-{1\over 2}} e_1+2^{-{1\over 2}} e_{n+1}$.
Consider $p_0$ in $A^{**}$ given by $(p_0)_\infty=0$ and $(p_0)_n=v_n\times v_n$ and $p(u)$ given by $p(u)_n=v_n\times v_n$ and $p(u)_\infty=u\times u$, where $u$ is a unit vector.
Then $p_0$ is an open projection, and $p(u)$ is closed if and only if $u=\lambda e_1$.
If $q$ is a projection such that $\|p(u)-q\|<1$, then $q_n$ has rank one, $1\leq n\leq\infty$.
It is then easy to see that $q$ is compact if and only if $q_n\to q_\infty$ in norm.

Of course $\alpha(p(u))=\alpha(p_0)=2$, and dist$(p(u),RC)=$ dist$(p_0,RC)=2^{-{1\over 2}}$, for any $u$.
dist$(p_0,CRC)=1$ by 5.1.
The determination of dist$(p(u),CRC)$ reduces to an elementary problem.
Clearly, for $q$ compact, $\|p(u)-q\|\geq \|p(u)_\infty-q_\infty\|$, $\limsup_n\|p(u)_n-q_\infty\|$.

If the $\limsup$ is $L$, we can modify the $q_n$'s so that $\|p(u)_n-q_n\|\leq L+\var$, $\forall n$.
(Actually, the ``$\var$'' is unnecessary.)
Also, if $q_\infty=w\times w$, $\|w\|=1$, then\newline
$\limsup d_a (p(u)_n,q_\infty)=\cos^{-1}|2^{-{1\over 2}} (e_1,w)|$, since $v_n\overset w\to\longrightarrow 2^{-{1\over 2}}e_1$.
Therefore\newline
$d_a(p(u),CRC)=\cos^{-1}\sup \{\min (|(u,w)|,2^{-{1\over 2}}|(e_1,w)|)\colon \|w\|=1\}$.
(Recall that dist$(\cdot,CRC)=\sin d_a (\cdot,CRC)$.)
Assume, as we may, that $(u,e_1)\geq 0$.
Then we can solve this maximin problem as follows:\ If $(u,e_1)\geq 2^{-{1\over 2}}$, let $w=e_1$.
If $(u,e_1) < 2^{-{1\over 2}}$, choose $w$ of the form $se_1+tu,s,t\geq 0$ so that $(u,w)=2^{-{1\over 2}}(e_1,w)$.

From the above we see that dist$(p(u),CRC) < 1$, $\forall u$, and dist$(p(u),CRC)=$ dist$(p(u),RC)$ if and only if $|(u,e_1)|\geq 2^{-{1\over 2}}$.
The largest value of dist$(p(u),CRC),({2\over 3})^{1\over 2}$, occurs when $(u,e_1)=0$.
The closure of $p(u)$ is given by $\overline{(p(u))_\infty}=(u\times u)\vee (e_1\times e_1)$, a rank two projection except when $u=\lambda e_1$.
It is easy to see that $\alpha(\overline{p(u))}=2$.
Thus if $|(u,e_1)| < 2^{-{1\over 2}}$, we have $p_0\leq p(u)\leq\overline{p(u)}$ and dist$(p_0,CRC) >$ dist$(p(u),CRC) >$ dist$(\overline{p(u)},CRC)$.

Let $\varphi_n$ be the pure state given by $\varphi_n(a)=(a_n v_n,v_n)$.
Then $\varphi_n\overset{w^*}\to\longrightarrow {1\over 2}\varphi$, where $\varphi(a)=(a_\infty e_1,e_1)$.
If $(u,e_1)=0$, then the support projection of $\varphi$ is orthogonal to $p(u)$, just as it is orthogonal to $p_0$.

It would seem that the study of dist$(p,CRC)$, for general $p$, is more complicated than the study of dist$(p,RC)$.
It would be interesting to know whether there is any natural hypothesis on $p$ (other than that $p$ be closed) which, together with $\alpha(p)<\infty$, implies dist$(p,CRC)<1$.

\bye